\documentclass[12pt,reqno]{amsart}
\usepackage{amssymb}

\usepackage{amscd}

\newtheorem{theorem}{Theorem}[section]
\newtheorem{corollary}[theorem]{Corollary}
\newtheorem{lemma}[theorem]{Lemma}
\newtheorem{proposition}[theorem]{Proposition}

\newtheorem{definition}[theorem]{Definition}
\newtheorem{remark}[theorem]{Remark}

\newtheorem{example}[theorem]{Example}
\numberwithin{equation}{section}
\begin{document}
\title{Perturbations of Dirac Operators}
\author{Igor Prokhorenkov}
\author{Ken Richardson}
\address{Department of Mathematics\\
Texas Christian University\\
Box 298900 \\
Fort Worth, Texas 76129}
\email{i.prokhorenkov@tcu.edu\\
k.richardson@tcu.edu}
\date{July, 2004}

\begin{abstract}
We study general conditions under which the computations of the index of a
perturbed Dirac operator $D_{s}=D+sZ$ localize to the singular set of the
bundle endomorphism $Z$ in the semi-classical limit $s\to \infty $. We show
how to use Witten's method to compute the index of $D$ by doing a
combinatorial computation involving local data at the nondegenerate singular
points of the operator $Z$. In particular, we provide examples of novel
deformations of the de Rham operator to establish new results relating the
Euler characteristic of a spin$^{c}$ manifold to maps between its even and
odd spinor bundles. The paper contains a list of the current literature on
the subject.
\end{abstract}


\maketitle

\section{Introduction}

About 20 years ago E. Witten \cite{Wit1} introduced a beautiful new approach
to proving Morse inequalities based on the deformation of the de Rham
complex. His ideas were fruitfully applied in many different situations
briefly outlined in the historical section at the end of the introduction.

\subsection{Example: Poincar\'{e}-Hopf theorem\label{PoincareSection}}

As a motivating example we sketch the Witten approach to the proof of the
Poincar\'{e}-Hopf theorem. Let $V$ be a smooth vector field on a smooth,
closed manifold $M$. A point $\overline{x}\in M$ is called a \textit{singular%
} point of $V$ if $V\left( \overline{x}\right) =0$ . In local coordinates $%
\left( x_{i}\right) $ near $\overline{x}$, we write 
\begin{equation*}
V\left( x\right) =\sum_{i}V^{i}\left( x\right) \partial _{i}.
\end{equation*}
A singular point $\overline{x}$ is called \textit{non-degenerate} if 
\begin{equation*}
\det \left( \frac{\partial V_{i}}{\partial x_{j}}\right) \left( \overline{x}%
\right) \neq 0.
\end{equation*}
The property of being non-degenerate does not depend on the system of
coordinates. All non-degenerate critical points are isolated, and there are
only finite number of them. The \textit{index} of $\overline{x}$ is defined
to be 
\begin{equation*}
\mathrm{ind}\left( \overline{x}\right) =\mathrm{sign}\det \left( \frac{%
\partial V_{i}}{\partial x_{j}}\right) \left( \overline{x}\right) .
\end{equation*}
Let $n^{\pm }$ denote the number of singular points with the index $\pm 1$.
Then the \textit{Poincar\'{e}-Hopf Theorem} asserts that 
\begin{equation}
\chi \left( M\right) =n^{+}-n^{-},  \label{PoincareHopf}
\end{equation}
where $\chi \left( M\right) $ is the Euler characteristic of $M.$

In order to prove the Poincar\'{e}-Hopf theorem, Witten\footnote{%
In \cite{Wit1} Witten mostly considered the case when $V$ is a Killing
vector field.} introduced the following one-parameter deformation of the
differential in the de Rham complex of $M$ depending on a real parameter $s$:

\begin{equation*}
d_{s}:=d+sV\lrcorner :\Omega ^{\bullet }\left( M\right) \rightarrow \Omega
^{\bullet }\left( M\right) ,
\end{equation*}
where $V\lrcorner $ denotes the interior product. Choose a Riemannian metric 
$g$ on $M$, and let $\left( d_{s}\right) ^{*}$ be the adjoint operator with
respect to the $L^{2}$ inner product on the space $\Omega ^{\bullet }\left(
M\right) $ of smooth forms. The corresponding deformation of the de Rham
operator $D=d+d^{*}$ is 
\begin{equation*}
D_{s}=d_{s}+d_{s}^{*}=d+d^{*}+s(V^{\flat }\wedge +V\lrcorner ),
\end{equation*}
where $V^{\flat }$ denotes the 1-form dual to $V$. The deformed de Rham
operator is a first order, elliptic differential operator that is in fact a
perturbed Dirac operator acting on the Clifford module of exterior forms on $%
M$ (see Section \ref{deRham}). We denote the restrictions of $D_{s}$ to even
or odd forms as $D_{s}^{\pm }:\Omega ^{\pm }\left( M\right) \rightarrow
\Omega ^{\mp }\left( M\right) $. The \textit{graded Witten Laplacian} is 
\begin{equation*}
\left( D_{s}^{2}\right) ^{\pm }:=D_{s}^{\mp }\circ D_{s}^{\pm }:\Omega ^{\pm
}\left( M\right) \rightarrow \Omega ^{\pm }\left( M\right) .
\end{equation*}
Since the index of the elliptic operator does not depend on the lower order
terms, for each $s$ we have 
\begin{eqnarray*}
\chi \left( M\right) &=&\mathrm{\mathrm{\mathrm{ind}}}\left( D_{s}\right)
=\dim \ker \left( D_{s}^{+}\right) -\dim \ker \left( D_{s}^{-}\right) \\
&=&\dim \ker \left( \left. \left( D_{s}\right) ^{2}\right| _{\Omega
^{+}\left( M\right) }\right) -\dim \ker \left( \left. \left( D_{s}\right)
^{2}\right| _{\Omega ^{-}\left( M\right) }\right) .
\end{eqnarray*}
An easy calculation shows that the Witten Laplacian $\left( D_{s}\right)
^{2} $ has the form 
\begin{equation*}
\left( D_{s}\right) ^{2}=\left( d_{s}+d_{s}^{*}\right) ^{2}=\left(
d+d^{*}\right) ^{2}+s^{2}|V|^{2}+sB,
\end{equation*}
where $B$ is a smooth bundle map. For very large $s$, the ``potential
energy''\ $s^{2}|V|^{2}$ becomes very large except in a small neighborhood
of the singular set of $V$. As explained in the original paper of Witten 
\cite{Wit1} and made rigorous in subsequent works by other researchers (see
next section for citations), the eigenforms of $\left( D_{s}\right) ^{2}$
concentrate near the singular points of $V$. There are asymptotic formulas
for the eigenvalues of the Witten Laplacian in terms of data at the singular
set of $V$. Local computations then complete the proof of the
Poincar\'{e}-Hopf theorem.

\subsection{\protect\smallskip The content of the paper}

The purpose of this paper is to study general conditions under which one can
use the method of Witten deformations to obtain an expression for the index
of a Dirac operator in terms of local quantities associated to a singular
set of a given bundle map. We now describe the setup of our paper. See
Section \ref{Preliminary} for a review of graded Clifford bundles and Dirac
operators.

Let $E=E^{+}\oplus E^{-}$ be a graded self-adjoint Clifford module over a
closed, smooth, Riemannian manifold $M$. Let $\Gamma \left( M,E\right) $
denote the space of smooth sections of $E$ and $D:\Gamma \left( M,E\right)
\rightarrow \Gamma \left( M,E\right) $ be the Dirac operator associated to a
Clifford module $E$.

Let $D^{\pm }:\Gamma \left( M,E^{\pm }\right) \rightarrow \Gamma \left(
M,E^{\mp }\right) $ denote the restrictions of the Dirac operator to smooth
even or odd sections. Observe that $D^{-}=\left( D^{+}\right) ^{\ast }$, the 
$L^{2}$-adjoint of $D^{+}$. Let $Z^{+}:\Gamma \left( M,E^{+}\right)
\rightarrow \Gamma \left( M,E^{-}\right) $ be a smooth bundle map, and we
let $Z^{-}$ denote the adjoint of $Z^{+}$. The operator $Z$ on $\Gamma
\left( M,E\right) $, defined by $Z\left( v^{+}+v^{-}\right)
=Z^{-}v^{-}+Z^{+}v^{+}$ for any $v^{+}\in E_{x}^{+}$ and $v^{-}\in E_{x}^{-}$
, is self-adjoint. A \textit{generalized Witten deformation} of $D$ is a
family $D_{s}$ of perturbed differential operators 
\begin{equation*}
D_{s}=\left( D+sZ\right) :\Gamma \left( M,E\right) \rightarrow \Gamma \left(
M,E\right) .
\end{equation*}
We define the operators $D_{s}^{\pm }$ by restricting in the obvious ways.
Our definition includes the known examples of Witten deformation as special
cases.

It is well-known \cite{A-S1} that the index $\mathrm{ind}\left( D^{+}\right) 
$ of $D^{+}$ depends only on the homotopy type of the principal symbol and
satisfies 
\begin{equation*}
\mathrm{ind}\left( D\right) =\dim \ker \left( \left. \left( D_{s}\right)
^{2}\right| _{\Gamma \left( M,E^{+}\right) }\right) -\dim \ker \left( \left.
\left( D_{s}\right) ^{2}\right| _{\Gamma \left( M,E^{-}\right) }\right) .
\end{equation*}
Thus, we need to study the operator 
\begin{equation*}
\left( D_{s}\right) ^{2}=D^{2}+s\left( ZD+DZ\right) +s^{2}Z^{2}.
\end{equation*}
The leading order behavior of the eigenvalues of this operator as $%
s\rightarrow \infty $ is determined by combinatorial data at the singular
set of the operator $Z$. This ``localization''\ allows one to compute $%
\mathrm{ind}\left( D\right) $ in terms of that data.

The results of this paper hold if the perturbation $Z$ is a \emph{proper
perturbation}. That is, it satisfies two conditions:

\begin{enumerate}
\item  \label{Zeroordercondition}$\left( D_{s}\right) ^{2}-\left( D\right)
^{2}=s\left( ZD+DZ\right) +s^{2}Z^{2}$ is a bundle map.

\item  \label{propPerturbCondition}Each singular point $\bar{x}$ of $Z$ is a 
\emph{proper singular point}; that is,

\begin{enumerate}
\item  $Z\left( \bar{x}\right) =0.$

\item  In local coordinates $x$ on a sufficiently small neighborhood $U$ of $%
\bar{x}$, there exists a constant $c>0$ such that for all $\alpha \in \Gamma
\left( U,E\right) $ and all $x\in U$, 
\begin{equation*}
\left\| Z\alpha \right\| _{x}\geq c\left| x-\bar{x}\right| ~\left\| \alpha
\right\| _{x},
\end{equation*}
where $\left\| \cdot \right\| _{x}$ is the pointwise norm on $E_{x}$.
\end{enumerate}
\end{enumerate}

We should note that in the important case when $Z$ is Clifford
multiplication by a vector field, Condition (\ref{Zeroordercondition}) is
not satisfied, and localization typically does not occur, (see Section \ref
{NoLocalization}).\footnote{%
If the vector field is a generator of an action by a one-parameter group of
isometries, the localization occurs when the Witten Laplacian is restricted
to each eigenspace of the Lie derivative associated to this generator. We
plan to treat this situation in a paper currently under preparation.} In the
proof of the Poincar\'{e}-Hopf index theorem, $Z=V^{\flat }\wedge
+V\lrcorner $, and Condition (\ref{propPerturbCondition}) reduces to the
requirement that the vector field has nondegenerate zeros.

In Section \ref{perturbDiracOpsSection} , we classify the possible gradings
of $E$ compatible with the existence of such $Z$. We also establish
necessary and sufficient conditions on the dimension of $E$ and on the form
of the operator $Z$ in order for (\ref{Zeroordercondition}) and (\ref
{propPerturbCondition}) to be satisfied.

In Section \ref{LocalizationSection}, we show that conditions (\ref
{Zeroordercondition}) and (\ref{propPerturbCondition}) imply that the
singular set of $Z$ consists of a finite number of non-degenerate zeros, and
as $s\rightarrow \infty $ the bounded spectrum of the Witten Laplacian
localizes to the singular set of $Z$. This means that the index of $D$ can
be computed by studying the zero spectrum of limiting ``model''\ operators,
which turn out to be harmonic oscillators. The main tool is the localization
theorem of M. Shubin \cite{Shu1}.

If a proper perturbation with no singularities exists, then the index is
zero. In particular, let $E^{\pm }\cong \left( \mathbb{S}\otimes W\right)
^{\pm }=\left( \mathbb{S}^{+}\otimes W^{\pm }\right) \oplus \left( \mathbb{S}%
^{-}\otimes W^{\mp }\right) $ be any graded, self-adjoint Clifford module
over an even dimensional, spin$^{c}$ manifold $M$ (all such Clifford modules
have this form; see Corollary \ref{GradedBundleClassification} in the
appendix). Then the index of the Dirac operator corresponding to this
Clifford module is zero if the bundles $W^{+}$ and $W^{-}$ are isomorphic
(see Corollary \ref{NoWIsomorphism}).

In Section \ref{localIndexSection}, we consider an elliptic operator of the
form $Q=\sum \left( A_{j}\partial _{j}+x_{j}B_{j}\right) $ on $\mathbb{C}%
^{m} $-valued functions on $\mathbb{R}^{n}$, such that each $A_{j}$ and $%
B_{k}$ is an $m\times m$ matrix and $\sum x_{k}B_{k}$ is a proper
perturbation. We prove that $Q$ is Fredholm and that continuous families of
such operators have the same index.

The main results of the paper are Theorem \ref{masterpiece0} and Theorem \ref
{Masterpiece}, which express the index of the Dirac operator $D$ in terms of
the local information at the singular points. We show in Theorem \ref
{masterpiece0} that the index of the Dirac operator $D$ is the sum of
indices of operators on vector-valued functions on $\mathbb{R}^{n}$ as in
Section \ref{localIndexSection}, where the coefficients $A_{j}$ and $B_{k}\ $%
depend only on local data at each singular point $\overline{x}$. Assuming
typical properties of $Z$ near each $\overline{x}$, the indices may be
computed more explicitly, as shown in Theorem \ref{Masterpiece}. We
emphasize strongly that all the information necessary to compute the index
of $D$ is contained in the set of matrices of first derivatives of $Z$ and
in the Clifford matrices taken at each singular point $\overline{x}$. Thus
this information is local in nature, and the answer is easily obtainable.

In Section \ref{ExamplesSection} we apply our results to the geometric Dirac
operators. In particular, we use Corollary \ref{EvenMasterpiece} to obtain
the Poincar\'{e}-Hopf theorem. Our proof also yields a new result, that the
Euler characteristic of an even dimensional, spin$^{c}$ manifold is zero if
and only if the even and odd spin$^{c}$ bundles are isomorphic (see
Corollary \ref{NoSpinIso}). In Theorem \ref{PinorHopf}, we show that the
Euler characteristic of a spin$^{c}$ manifold is the sum of the indices of
zeros of a possibly singular section of the conformal pin bundle over the
manifold. Thus, the Euler characteristic is zero if and only if the odd pin
bundle $\mathrm{Pin}^{-}\left( T^{*}M\right) $ has a global section. In
Section \ref{submanifoldSection}, we use our results to show that if $M$ is
a submanifold of odd codimension in a manifold endowed with a graded
Clifford module, then the index of the Dirac operator associated to the
restriction of this Clifford bundle to M is zero.

\subsection{Review of literature\label{Literature}}

Witten deformation was first introduced in \cite{Wit1}, where the author
sketched a beautiful proof of the Morse inequalities by deforming the de
Rham complex (see also \cite{Bott}, \cite{C-F-K-Si}, \cite{Hel1}, \cite{Hen}
, \cite{Shu1}, \cite{Si},\cite{Roe}). In addition, using the ideas of
quantum field theory (supersymmetry and instantons) Witten explained how to
derive analytically the entire Morse complex. In \cite{Hel-Sj1} B. Helffer
and J. Sj\"{o}strand put Witten's analysis on a rigorous footing. See also
current introductions to the subject in \cite{Bur}, \cite{Bur1}, and \cite
{C-P-R}. A recent discussion of connections between the spectral theory and
semiclassical analysis of the Witten Laplacian and the notion of
hypoellipticity is given in the lecture notes by B. Helffer \cite{Hel4}. In
addition Witten suggested a way to use his method to prove the
Poincar\'{e}-Hopf theorem. Rigorous treatments of his ideas in this
direction are contained in \cite{ES-W}, \cite{Shu2}, \cite{Yu}, and \cite
{Zh2}.

Witten obtained holomorphic Morse inequalities by the same method in \cite
{Wit2}. The asymptotic holomorphic Morse inequalities, were proved by J.-P.
Demailly and J.-M Bismut (see \cite{D1}, \cite{D2}, \cite{Siu}, \cite{Bis2}
). Equivariant holomorphic Morse inequalities were investigated by V.
Mathai, S. Wu, and W.Zhang in a series of papers \cite{Mat-Wu}, \cite{Wu},
and \cite{Wu-Zh}.

In \cite{Bis1} J.-M. Bismut modified the Witten deformation technique and
combined it with intricate and deep ideas of probability theory to produce a
new proof of the degenerate Morse-Bott inequalities. A more accessible
proof, based on the adiabatic technique of Mazzeo-Melrose and Forman (\cite
{Maz-Me} , \cite{Fo1}) instead of probability considerations, was given by
I. Prokhorenkov in \cite{Pro} ; see also \cite{Bra-Far1}, \cite{Hel-Sj2},
and \cite{Ji} for different proofs and generalizations.

A. V. Pazhitnov \cite{Pa} used the method of Witten to prove some of the
Morse-Novikov inequalities --- that is, when the gradient of Morse function
is replaced in the deformation by a closed, nondegenerate one-form \cite
{No-Shu}. Novikov inequalities for vector fields were established by M. A.
Shubin in an influential paper \cite{Shu2}. Shubin's results were extended
by M. Braverman and M. Farber \cite{Bra-Far1} to the case when the one-form
(or corresponding vector field) has non-isolated zeros, and by H. Feng and
E. Guo to the case of more general vector fields in \cite{F-G1} and \cite
{F-G2}. The equivariant Morse-Novikov inequalities were proved in \cite
{Bra-Far2} (see also \cite{Wai}). The Novikov-Witten complex was constructed
by D. Burghelea and S. Haller in \cite{B-H}. Holomorphic Morse inequalities
for K\"{a}hler manifolds in the presence of a holomorphic circle action and
their applications to symplectic reduction were studied by M. Braverman in 
\cite{Bra2}. M. Calaza \cite{Ca} proved a generalization of Morse
inequalities for an orbit space.

J. Alvarez L\'{o}pez \cite{A-L} used the method of Witten to prove Morse
inequalities for the invariant cohomology of the space of orbits of a
pseudogroup of isometries with applications to the basic cohomology of
Riemann foliations. V. Belfi, E. Park, and K. Richardson \cite{B-P-R} used
the Witten deformation of the basic Laplacian to prove an analog of the
Poincar\'{e}-Hopf index theorem for Riemannian foliations. It should be
noted that there are no known proofs of these foliation results by methods
that do not use the Witten deformation technique.

A combinatorial analog of Witten deformation was introduced by R. Forman in 
\cite{Fo2} and \cite{Fo3}. Forman's ideas were extended by V. Mathai and S.
Yates in \cite{Mat-Ya} to the case of infinite cell complexes, thus
obtaining discrete analogs of Morse $L^{2}$ inequalities.

Witten deformation proved to be very productive in studying relations
between analytic and combinatorial torsions. Deep results in these direction
were obtained in \cite{Bis-Zh1}, \cite{Bis-Zh2}, \cite{Bis-G}, \cite{Bra1}, 
\cite{Bra4}, \cite{B-F-K1}, \cite{B-F-K2}, and \cite{B-F-K-M}.

In \cite{Shu1}, \cite{Shu3}, and \cite{Mat-Shu} M. A. Shubin and V. Mathai
used the Witten deformation $\left( d_{s}+d_{s}^{*}\right) ^{2}$ of the
Laplacian to study $L^{2}$ Morse inequalities for regular covering manifolds
and flat Hilbert bundles over compact manifolds. A nice survey of these
results can be found in \cite{Mat}. The first paper \cite{Shu1} also
contains a very useful analysis of model operators appearing as limits of
localizing deformed operators. A recent preprint \cite{K-M-S} of Y.
Kordyukov, V. Mathai, and M. Shubin extends the analysis to the case of
projectively invariant elliptic operators with invariant Morse type
potentials on covering spaces of compact manifolds.

Some interesting applications to manifolds with negative curvature are
discussed in \cite{El-R} and \cite{Li} . Here again, other approaches are
known to work.

Applications of Witten analysis to statistical physics are discussed in \cite
{Hel2}, \cite{Hel3}, \cite{Jo}, \cite{Sj}, and \cite{Wang}.

An approach to the Hodge theory on topologically tame non-compact manifolds
by means of a Witten Laplacian with a potential rapidly increasing at
infinity was suggested by E. Bueler and I. Prokhorenkov in \cite{Bue}, \cite
{Pro}, and \cite{Bue-Pro}. Further results about Witten Laplacians on
non-compact manifolds were obtained in \cite{Ah-St}, \cite{Far-Sh}, \cite
{G-W1}, \cite{G-W2}, and \cite{G-W3}.

W. Zhang \cite{Zh1} used the following modification of Witten deformation 
\begin{equation*}
D_{V}=\frac{1}{2}\left( c(V)(d+d^{*})-(d+d^{*})c(V)\right) ,
\end{equation*}
where $c(V)$ denotes Clifford multiplication by a nowhere zero vector field $%
V$, to study the Kervaire semicharacteristic of odd-dimensional compact
manifolds. For a recent survey of these and other applications of Witten
deformations, see \cite{Zh2}.

\section{Perturbing Dirac operators\label{perturbDiracOpsSection}}

\subsection{Preliminaries and Notational Conventions\label%
{preliminaryperturb}}

Throughout this paper, the manifold $M$ is always assumed to be a smooth,
closed, oriented Riemannian manifold of dimension $n$, and $E$ $=E^{+}\oplus
E^{-}$ is assumed to be a graded, self-adjoint, Hermitian Clifford module
over $M$. If $M$ is spin$^{c}$, then $\mathbb{S}$ always denotes a complex
spinor bundle over $M$, a particular example of such a bundle $E$. We define 
$D:\Gamma \left( M,E\right) \rightarrow \Gamma \left( M,E\right) $ to be the
corresponding Dirac operator, and let $D^{\pm }:\Gamma \left( M,E^{\pm
}\right) \rightarrow \Gamma \left( M,E^{\mp }\right) $ denote the
restrictions of $D$ to smooth even and odd sections. The operator $D^{-}$ is
the adjoint of $D^{+}$ with respect to the $L^{2}$-metric on $\Gamma \left(
M,E\right) $\textbf{\ }defined by the Riemannian metric on $M$\ and the
Hermitian metric on $E$.

In the following, we denote Clifford multiplication by an element $v\in
T_{x}M$ on the fiber $E_{x}$ by $c\left( v\right) $. Clifford multiplication
by cotangent vectors will use the same notation: $c\left( \alpha \right)
:=c\left( \alpha ^{\#}\right) $, where $T_{x}^{*}M\overset{\#}{\rightarrow }%
T_{x}M$ is the metric isomorphism.

The natural grading on $E$ is induced by the action of the\textbf{\ }%
chirality operator $\gamma $. Recall that if $e_{1},...,e_{n}$ is an
oriented orthonormal basis of $T_{x}M$, then the chirality operator is
multiplication by 
\begin{equation*}
\gamma =i^{k}c(e_{1})...c(e_{n})\in \mathrm{End}\left( E_{x}\right) \mathbf{,%
}
\end{equation*}
where $k=n/2$ if $n$ is even and $k=\left( n+1\right) /2$ if $n$ is odd. In
this paper we study the other possible gradings as well. See the appendix
(Section \ref{Preliminary}) for more information.

Let $Z^{+}\in \Gamma \left( M,\mathrm{Hom}\left( E^{+},E^{-}\right) \right) $
be a smooth bundle map, and let $Z^{-}$ denote the adjoint of $Z^{+}$. The
operator $Z$ on $\Gamma \left( M,E\right) $,$\ $defined by $Z\left(
v^{+}+v^{-}\right) =Z^{-}v^{-}+Z^{+}v^{+}$ for any $v^{+}\in E_{x}^{+}$ and $%
v^{-}\in E_{x}^{-},$ is self-adjoint. Let $D_{s}$ denote the perturbed Dirac
operator 
\begin{equation}
D_{s}=\left( D+sZ\right) :\Gamma \left( M,E\right) \rightarrow \Gamma \left(
M,E\right) ,  \label{perturbed operator}
\end{equation}
and define the operators $D_{s}^{\pm }$ by restricting in the obvious ways.

\subsection{ Nonexistence of perturbations compatible with the natural
grading}

For any differential operator $L$, let $\sigma _{L}$ denote its principal
symbol. We will start with the result that holds for general first-order
operators.

\begin{lemma}
Let $L:\Gamma \left( M,E\right) \rightarrow \Gamma \left( M,E\right) $ be a
smooth, first-order differential operator, and let $Z:\Gamma \left(
M,E\right) \rightarrow \Gamma \left( M,E\right) $ be a bundle endomorphism.
Then the operator $LZ+ZL$ is a bundle map if and only if $Z\circ \sigma
_{L}\left( x,\xi \right) +\sigma _{L}\left( x,\xi \right) \mathbf{\circ }Z\
=0$ on $E_{x}$ for every $x\in M$, $\xi \in T_{x}^{*}M$.
\end{lemma}

\begin{proof}
The differential operator $ZL+LZ$ is zeroth order if and only if it commutes
with multiplication $m_{f}$ by any smooth function $f$ on $M$. We calculate
the commutator 
\begin{eqnarray*}
\left[ ZL+LZ,m_{f}\right] &=&Z\left[ L,m_{f}\right] +\left[ L,m_{f}\right]
Z\,\,\text{since }Z\text{ is zeroth order} \\
&=&\,i\left( Z\circ \sigma _{L}\left( df\right) +\sigma \left( L\right)
\left( df\right) \mathbf{\circ }Z\right) ,
\end{eqnarray*}
where for any $1$-form $\alpha $ on $M$, $\sigma _{L}\left( \alpha \right) $
is the bundle endomorphism defined by $\left. \sigma _{L}\left( \alpha
\right) \right| _{x}=\sigma _{L}\left( x,\alpha _{x}\right) $.
\end{proof}

Since 
\begin{equation}
\left( D_{s}\right) ^{2}-D^{2}=s\left( ZD+DZ\right) +s^{2}Z^{2},
\label{DsSquared}
\end{equation}
we have the following corollary:

\begin{corollary}
\label{zerocondition}For any $s\neq 0$ the operator $\left( D_{s}\right)
^{2}-D^{2}$ is zeroth order if and only if $Z\circ \sigma _{D}\left( x,\xi
\right) +\sigma _{D}\left( x,\xi \right) \mathbf{\circ }Z\ =0$ on $E_{x}$
for every $x\in M$, $\xi \in T_{x}^{*}M$.
\end{corollary}

\begin{remark}
Corollary \ref{zerocondition} is true in even greater generality, such as
when $D$ is a first-order, classical pseudodifferential operator. However,
in this paper we consider only differential operators.
\end{remark}

A bundle endomorphism $Z$ satisfying the condition in Corollary \ref
{zerocondition} does not always exist. In particular, the following result
applies to the spin$^{c}$\ Dirac operator, whose principal symbol is $%
ic\left( \xi \right) $ (for $\xi \in T^{*}M$).

\begin{proposition}
\label{NoZonS}Let $V$ be an even-dimensional, oriented, Euclidean vector
space. Let $\mathbb{S}=\mathbb{S}^{+}\oplus \mathbb{S}^{-}$ denote the
associated space of complex spinors. There does not exist a linear map $Z:%
\mathbb{S}\rightarrow \mathbb{S}$ such that $Z\circ c\left( v\right)
+c\left( v\right) \mathbf{\circ }Z\ =0$ for every $v\in V$ and such that $Z$
restricts to a map $Z:\mathbb{S}^{+}\rightarrow \mathbb{S}^{-}$.
\end{proposition}

\begin{proof}
Any endomorphism of $\mathbb{S}$ can be written as Clifford multiplication
by an element of $\mathbb{C}\mathrm{l}\left( V\right) $, so that the result
is equivalent to the statement that no element of $\mathbb{C}\mathrm{l}%
\left( V\right) $ anticommutes with every vector. This is a consequence of
the elementary fact that if $\alpha \in \mathbb{C}\mathrm{l}\left( V\right) $
anticommutes with every vector, then $\alpha $ is a complex multiple of the
chirality operator. Since this element maps $\mathbb{S}^{+}$ to itself, the
result follows.
\end{proof}

The corollary below is a generalization.

\begin{corollary}
\label{NoZonSxE}Let $V$ be an even-dimensional, oriented, Euclidean vector
space. Let $\mathbb{S}=\mathbb{S}^{+}\oplus \mathbb{S}^{-}$ denote the
associated space of complex spinors. Let $W$ be another vector space, and
consider the Clifford action $\widetilde{c}\left( v\right) =c\left( v\right)
\otimes \mathbf{1}$ on $\mathbb{S}\otimes W$. There does not exist a linear
map $Z:\mathbb{S}\otimes W\rightarrow \mathbb{S}\otimes W$ such that $Z\circ 
\widetilde{c}\left( v\right) +\widetilde{c}\left( v\right) \mathbf{\ \circ }%
Z\ =0$ for every $v\in V$ and such that $Z$ restricts to a map $Z:\mathbb{S}
^{+}\otimes W\rightarrow \mathbb{S}^{-}\otimes W$.
\end{corollary}

\begin{proof}
Endow $W$ with a Euclidean (or Hermitian) metric . Choose an orthonormal
basis $\left\{ w_{1},...,w_{k}\right\} $ of $W$. Write 
\begin{equation*}
Z(v\otimes w_{i})=\sum_{j}Z_{ij}(v)\otimes w_{j}\,;
\end{equation*}
the corollary follows from the proposition applied to each linear operator $%
Z_{ij}$ separately.
\end{proof}

\begin{example}
The Dolbeault and signature operators do not have such perturbations,
because in both cases the Clifford action has the form of Corollary \ref
{NoZonSxE}.
\end{example}

\subsection{\protect\smallskip Admissible Perturbations}

The following results determine the precise form of perturbations $Z$
satisfying the condition in Corollary \ref{zerocondition} if $D$ is the
Dirac operator associated to a Clifford bundle over a spin$^{c}$ manifold.

\begin{proposition}
\label{SpinC}Let $M$ be even-dimensional and spin$^{c}$. Let $E^{\pm }\cong
\left( \mathbb{S}\otimes W\right) ^{\pm }=\left( \mathbb{S}^{+}\otimes
W^{\pm }\right) \oplus \left( \mathbb{S}^{-}\otimes W^{\mp }\right) $ be any
Clifford module over $M$ (see Corollary \ref{GradedBundleClassification} in
the appendix). Suppose that there is a bundle endomorphism $Z^{+}:\Gamma
\left( M,\left( \mathbb{S}\otimes W\right) ^{+}\right) \rightarrow \Gamma
\left( M,\left( \mathbb{S}\otimes W\right) ^{-}\right) $ such that the
self-adjoint operator $Z=\left( Z^{+},\left( Z^{+}\right) ^{*}\right) :%
\mathbb{S}\otimes W\rightarrow \mathbb{S}\otimes W$ anticommutes with
Clifford multiplication by vectors. Then $Z$ has the form $Z=\gamma \otimes
\phi $, where $\gamma $ is the chirality operator on $\mathbb{S}$ and where $%
\phi ^{+}:W^{+}\rightarrow W^{-}$ is a bundle map with $\phi =\left( \phi
^{+},\left( \phi ^{+}\right) ^{*}\right) $. Conversely, any bundle
endomorphism of that form anticommutes with Clifford multiplication by
vectors.
\end{proposition}

\begin{proof}
The action of $Z$ on $\mathbb{S}\otimes W$ has the following local form. For
a local orthonormal basis $\left\{ b_{1},...,b_{k}\right\} $ of $W$and any $%
\alpha \in \Gamma \left( M,\mathbb{S}\right) $, 
\begin{equation*}
Z(\alpha \otimes b_{i})=\sum_{j}Z_{ij}(\alpha )\otimes b_{j},
\end{equation*}
where each operator $Z_{ij}$ must anticommute with Clifford multiplication
by vectors. Thus $Z_{ij}=c_{ij}\,\gamma (\alpha )$ for some complex scalar $%
c_{ij}$ and 
\begin{equation*}
Z(\alpha \otimes b_{i})=\sum_{j}c_{ij}\,\gamma (\alpha )\otimes b_{j}=\gamma
(s)\otimes \sum_{j}c_{ij}\,b_{j}\,.
\end{equation*}
The operator $\gamma $ restricts to $\mathbf{1}$ on $\mathbb{S}^{+}$ and $%
\mathbf{-1}$ on $\mathbb{S}^{-}$. The result follows from the hypothesis and
the equation above.
\end{proof}

\begin{remark}
The case when $W^{-}$ is zero-dimensional is reflected in Corollary \ref
{NoZonSxE}, where the only possible bundle map $\phi $ is the zero map.
\end{remark}

\begin{remark}
The restriction of such a $\phi $ to a fiber is invertible if and only if $%
\mathrm{\dim }W^{+}=\mathrm{\dim }W^{-}$.
\end{remark}

\begin{proposition}
\label{SpinCOdd}Let $M$ be odd dimensional and spin$^{c}$. Fix the
representation $c^{+}$ of $\mathbb{C}\mathrm{l}\left( TM\right) $. Let $
E\cong \mathbb{S}\otimes \left( W^{\prime }\oplus W^{\prime }\right) $ be as
in Corollary \ref{GradedBundleClassification}. Suppose that there exists a
self-adjoint endomorphism $Z=\left( Z^{+},\left( Z^{+}\right) ^{*}\right)
:E\rightarrow E$ that anticommutes with the Clifford multiplication $\left(
c^{+}\left( v\right) \otimes \mathbf{1},\,c^{-}\left( v\right) \otimes 
\mathbf{1}\right) $ by all vectors $v\in TM$. Then $Z$ has the form $Z=%
\mathbf{1}\otimes \left( 
\begin{array}{ll}
0 & \phi \\ 
-\phi & 0
\end{array}
\right) $, where $\phi :W^{\prime }\rightarrow W^{\prime }$ is a
skew-adjoint bundle map. Conversely, any bundle endomorphism of that form is
self-adjoint and anticommutes with Clifford multiplication by vectors.
\end{proposition}

\begin{proof}
Since $\mathbb{C}\mathrm{l}^{+}\cong \mathrm{End}\left( \mathbb{S}\right) $
in odd dimensions, no nonzero element of $\mathrm{End}\left( \mathbb{S}%
\right) $ anticommutes with Clifford multiplication $c^{+}$ by all vectors.
An endomorphism $Z$ that anticommutes with 
\begin{eqnarray*}
c^{E}\left( v\right) &=&\left( c^{+}\left( v\right) \otimes \mathbf{1}
,c^{-}\left( v\right) \otimes \mathbf{1}\right) \\
&=&c^{+}\left( v\right) \otimes \left( 
\begin{array}{cc}
\mathbf{1} & 0 \\ 
0 & -\mathbf{1}
\end{array}
\right)
\end{eqnarray*}
for all vectors $v$ must have the following properties. First, $Z:\mathbb{S}%
\otimes W\rightarrow \mathbb{S}\otimes W$ must commute with all maps of the
form $A\otimes \mathbf{I}$ with $A\in \mathrm{End}\left( \mathbb{S}\right) $
, because any such map can be realized as multiplication by an element in $%
\mathbb{C}\mathrm{l}^{+}$ in odd dimensions. It follows easily that $Z$ can
be expressed in the form $Z=\mathbf{1}\otimes Z^{\prime }$, where $Z^{\prime
}=\left( 
\begin{array}{cc}
Z_{1} & Z_{2} \\ 
Z_{3} & Z_{4}
\end{array}
\right) $ is an endomorphism of $W^{\prime }\oplus W^{\prime }$. Because $%
Z^{\prime }$ maps $\mathrm{span}\left\{ \left( w,\pm w\right) \in W^{\prime
}\oplus W^{\prime }\,|\,w\in W^{\prime }\right\} $ to $\mathrm{span}\left\{
\left( w,\mp w\right) \in W^{\prime }\oplus W^{\prime }\,|\,w\in W^{\prime
}\right\} $, $Z_{1}\pm Z_{2}=\mp Z_{3}-Z_{4}$. These equations imply that $%
Z_{4}=-Z_{1}$ and $Z_{3}=-Z_{2}$. Because $Z$ anticommutes with $c^{E}\left(
v\right) $, $Z_{1}=Z_{4}=0$. Self-adjointness of $Z$ implies further that $%
Z_{2}$ is skew-adjoint, and the result follows.
\end{proof}

\begin{remark}
If $M$ is not spin$^{c}$, Propositions \ref{SpinC} and \ref{SpinCOdd} above
remain true locally.
\end{remark}

\subsection{ Proper perturbations of Dirac Operators\label%
{properPerturbSection}}

In this section, we state the nondegeneracy conditions on the perturbation.
In previously studied types of perturbations, our conditions are equivalent
to those required by others (see the introduction).

\begin{definition}
Let $Z:E\rightarrow E$ be a smooth bundle map. We say that $\overline{x}\in
M $ is a \textbf{proper singular point} of $Z$ if on a sufficiently small
neighborhood $U$ of $\overline{x}$ we have

\begin{enumerate}
\item  $Z|_{\overline{x}}=0$, and

\item  in local coordinates $x$ on $U$, there exists a constant $c>0$ such
that for all $\alpha \in \Gamma \left( U,E\right) $, 
\begin{equation*}
\left\| Z\alpha \right\| _{x}\geq c\left| x-\overline{x}\right| \left\|
\alpha \right\| _{x},
\end{equation*}
where $\left\| \,\cdot \,\right\| _{x}$ is the pointwise norm on $E_{x}$.
\end{enumerate}
\end{definition}

\begin{lemma}
A point $\overline{x}\in M$ is a proper singular point of $Z$ if and only
if, in local coordinates $x$ on $U$, there exist invertible bundle maps $%
Z_{j}$ for $1\leq j\leq n=\dim M$ over $U$ such that $Z=\sum_{j}\left( x-%
\overline{x}\right) _{j}Z_{j}$ on $U$, and $Z$ is invertible over $%
U\setminus \left\{ \overline{x}\right\} $.
\end{lemma}

\begin{proof}
Since $Z$ is smooth and vanishes at $\overline{x}$, $Z=\sum_{j}\left( x-%
\overline{x}\right) _{j}Z_{j}$ for some bundle maps $Z_{j}$. The inequality
in the definition above is equivalent to 
\begin{equation*}
\left\| \sum_{j}\sigma _{j}Z_{j}\alpha \right\| _{x}^{2}\geq c^{2}
\end{equation*}
for every $\sigma \in S^{n-1},$ $\alpha \in \Gamma \left( \overline{U}%
,E\right) ,$ $x\in U$ such that $\left\| \alpha \right\| _{x}=1$. Since the
left hand side of the inequality is a continuous function of $\sigma $ and $%
\alpha $ over the compact set $S^{n-1}\times \left\{ \alpha \in \Gamma
\left( \overline{U},E\right) \,|\,\left\| \alpha \right\| _{x}=1\text{ for
all }x\in U\right\} $, its infimum is attained. It follows that $Z$ is
invertible away from $\overline{x}$ if and only if the inequality holds.
\end{proof}

\begin{definition}
\label{properPerturbationDef}Let $D^{\pm }:\Gamma \left( M,E^{\pm }\right)
\rightarrow \Gamma \left( M,E^{\mp }\right) $ be the Dirac operator
associated to a bundle of graded Clifford modules. Let $D_{s}=D+sZ$ for $%
s\in \mathbb{R}$, where $Z=\left( Z^{+},\left( Z^{+}\right) ^{\ast }\right)
\in \Gamma \left( M,\mathrm{End}\left( E^{+},E^{-}\right) \right) $. We say
that $Z$ is a \textbf{proper perturbation of }$D$ if

\begin{enumerate}
\item  $\left( D_{s}\right) ^{2}-D^{2}$ \ is a zeroth order operator.

\item  All singular points of $Z$ are proper.
\end{enumerate}
\end{definition}

\begin{remark}
If $M$ is compact, then the number of singular points of a proper
perturbation is finite.
\end{remark}

The following lemma will be used to quantify ranks of vector bundles on
which the proper perturbations act.

\begin{lemma}
\label{dimensionLemma}There exists a linear map $L:\mathbb{R}^{k}\rightarrow
M_{r}\left( \mathbb{C}\right) $ that satisfies $L\left( x\right) \in \mathrm{%
Gl}\left( r,\mathbb{C}\right) $ for $x\neq 0$ if and only if $%
r=m2^{\left\lfloor \frac{k-1}{2}\right\rfloor }$ for some positive integer $%
m $.
\end{lemma}

\begin{proof}
Since $L\left( x\right) $ could be Clifford multiplication by the vector $x$
on the $r$-dimensional vector space $\mathbb{S}\otimes \mathbb{C}^{m}$,
ranks of the form $r=m2^{\left\lfloor \frac{k}{2}\right\rfloor }$ are
realizable. If $k$ is even, the image of the restriction of Clifford
multiplication to $\mathbb{S}^{+}\otimes \mathbb{C}^{m}$ has rank $r=m2^{%
\frac{k}{2}-1}$. Hence for all positive integers $k$ and $m$ , there exist
linear maps $L:\mathbb{R}^{k}\rightarrow M_{r}\left( \mathbb{C}\right) $
with $r=m2^{\left\lfloor \frac{k-1}{2}\right\rfloor }$ and $L\left( x\right)
\in \mathrm{Gl}\left( r,\mathbb{C}\right) $ for $x\neq 0$.

\footnote{%
This part of the argument is provided in \cite{Re} and \cite{Re2}.}Next,
suppose that a linear map $L:\mathbb{R}^{k}\rightarrow M_{r}\left( \mathbb{C}%
\right) $ satisfies $L\left( x\right) \in \mathrm{Gl}\left( r,\mathbb{C}%
\right) $ for $x\neq 0$, for some positive integers $k$ and $r$. Such a map
restricts to a map $L:S^{k-1}\rightarrow $ $\mathrm{Gl}\left( r,\mathbb{C}%
\right) $ with $L\left( -x\right) =-L\left( x\right) $. Consider the vector
bundles $r\mathbf{T}=\mathbb{R}P^{k-1}\times \mathbb{C}^{r}$ and $rL_{k}$ ($%
r $ times the complexification of the canonical line bundle) over the
projective space $\mathbb{R}P^{k-1}$. Note that 
\begin{eqnarray*}
\mathbb{R}P^{k-1}\times \mathbb{C}^{r} &=&S^{k-1}\times \mathbb{C}%
^{r}\diagup \left( x,y\right) \sim \left( -x,y\right) ,\text{ and} \\
rL_{k} &=&S^{k-1}\times \mathbb{C}^{r}\diagup \left( x,y\right) \sim \left(
-x,-y\right) .
\end{eqnarray*}
The map $f:S^{k-1}\times \mathbb{C}^{r}\rightarrow S^{k-1}\times \mathbb{C}%
^{r}$ defined by $f\left( x,y\right) =\left( x,L\left( x\right) y\right) $
induces an isomorphism between $rL_{k}$ and $r\mathbf{T}$. Thus, the virtual
bundle $r\left( L_{k}-\mathbf{T}\right) $ represents the zero element in the
reduced complex K-group $\widetilde{K}\left( \mathbb{R}P^{k-1}\right) $.
Since 
\begin{equation*}
\widetilde{K}\left( \mathbb{R}P^{k-1}\right) \cong \mathbb{Z}%
_{2^{\left\lfloor \frac{k-1}{2}\right\rfloor }}
\end{equation*}
with generator $L_{k}-\mathbf{T}$ (see \cite{Ad}), we must have $%
\displaystyle r=m2^{\left\lfloor \frac{k-1}{2}\right\rfloor }$ for some
positive integer $m$.
\end{proof}

The following two theorems give necessary and sufficient conditions on the
bundle $E$ and the bundle map $Z$ in order that $Z$ be a proper perturbation.

\begin{theorem}
\label{EvenProperPerturbations}Suppose that the dimension $n$ of $M$ is
even. Let $Z$ be a proper perturbation of $D$ on $\Gamma \left( M,E\right) $%
, with notation as in Definition \ref{properPerturbationDef}. Let $\left\{ 
\overline{x}_{1},\overline{x}_{2},...,\overline{x}_{N}\right\} $ be the
singular points of $Z$, and let $U_{l}$ be the neighborhood of $\overline{x}%
_{l}$ as in the definition. Over each $U_{l}$, choose a local spin$^{c}$
bundle and isomorphism $E\cong \mathbb{S}\otimes W$ as in Proposition \ref
{SpinC}. Then the following conditions must be satisfied.

\begin{enumerate}
\item  The local bundles $W^{+}$ and $W^{-}$ from Proposition \ref{SpinC}
must have the same dimension, which implies that the rank of $E$ must be a
multiple of $2^{\frac{n}{2}+1}$. If the set of singular points is nonempty,
then the rank must be of the form $m2^{n}$, where $m$ is a positive integer.

\item  Near each singular point $\overline{x}$, the bundle map $Z$ has the
form $Z=\sum_{j}\left( x-\overline{x}\right) _{j}\gamma \otimes \phi _{j}$,
with notation as in Definition \ref{properPerturbationDef} and Proposition 
\ref{SpinC}, where each $\phi _{j}^{+}:W^{+}\rightarrow W^{-}$ is a locally
defined bundle isomorphism with $\phi _{j}=\left( \phi _{j}^{+},\left( \phi
_{j}^{+}\right) ^{*}\right) $.
\end{enumerate}

Conversely, every graded, self-adjoint bundle map $Z$ that has a finite set
of singular points and satisfies the two conditions above is a proper
perturbation.
\end{theorem}

\begin{proof}
The first part of the first condition and the second condition follow
directly from Proposition \ref{SpinC}. To prove the first condition in the
case where the set of singular points is nonempty, suppose that $\overline{x}
$ is a singular point of $Z$, and choose coordinates centered at $\overline{x%
}$. Identify $W_{\overline{x}}^{+}\cong W_{\overline{x}}^{-}\cong \mathbb{C}
^{d}$ for the appropriate positive integer $d$. Define $L\left( x\right)
=\sum_{j=1}^{n}x_{j}\phi _{j}^{+}\left( \overline{x}\right) :\mathbb{C}%
^{d}\rightarrow \mathbb{C}^{d}$; the hypotheses imply that $L\left( x\right) 
$ is invertible for each $x\neq 0$. Thus, 
\begin{equation*}
d=m2^{\left\lfloor \frac{n-1}{2}\right\rfloor }=m2^{\frac{n}{2}-1}
\end{equation*}
for some positive integer $m$, by Lemma \ref{dimensionLemma}. Since the
dimension of $\mathbb{S}$ is $2^{\frac{n}{2}}$, the statement follows. The
converse is clear.
\end{proof}

\begin{remark}
The given minimal rank $2^{\frac{n}{2}+1}$ is sharp, since the example given
in Proposition \ref{Cobordism} has precisely that rank.
\end{remark}

\begin{theorem}
\label{OddProperPerturbations}Suppose that the dimension $n$ of $M$ is odd.
Let $Z$ be a proper perturbation of $D$ on $\Gamma \left( M,E\right) $, with
notation as in Definition \ref{properPerturbationDef}. Let $\left\{ 
\overline{x}_{1},\overline{x}_{2},...,\overline{x}_{N}\right\} $ be the
singular points of $Z$, and let $U_{l}$ be the neighborhood of $\overline{x}%
_{l}$ as in the definition. Over each $U_{l}$, choose a local spin$^{c}$
bundle and isomorphism $E\cong \mathbb{S}\otimes \left( W^{\prime }\oplus
W^{\prime }\right) $ as in Proposition \ref{SpinCOdd}. Then the following
condition must be satisfied.

\begin{enumerate}
\item  If the set of singular points is empty, there is no further
restriction on the rank of $E$; that is, it need only be a multiple of $2^{%
\frac{\left( n+1\right) }{2}}$. If the set of singular points is nonempty,
then the rank of $E$ must have the form $m2^{n}$, where $m$ is a positive
integer.

\item  Near each singular point, the bundle map $Z$ has the form $%
Z=\sum_{j}x_{j}\mathbf{1}\otimes \left( 
\begin{array}{ll}
0 & \phi _{j} \\ 
-\phi _{j} & 0
\end{array}
\right) $ as in Proposition \ref{SpinCOdd} and in Definition \ref
{properPerturbationDef}, where each $\phi _{j}:W^{\prime }\rightarrow
W^{\prime }$ is a locally defined, skew-adjoint bundle isomorphism.
\end{enumerate}

Conversely, every graded, self-adjoint bundle map $Z$ that has a finite set
of singular points and satisfies the two conditions above is a proper
perturbation.
\end{theorem}

\begin{proof}
The second condition follows directly from Proposition \ref{SpinCOdd}. To
prove the first condition in the case where the set of singular points is
nonempty, choose local coordinates $x$ centered at a singular point $%
\overline{x}$. Identify $W_{\overline{x}}^{\prime }\cong \mathbb{C}^{d}$ for
the appropriate positive integer $d$. Define $L\left( x\right)
=\sum_{j=1}^{n}x_{j}\phi _{j}\left( \overline{x}\right) :\mathbb{C}%
^{d}\rightarrow \mathbb{C}^{d}$; the hypotheses imply that $L\left( x\right) 
$ is invertible for each $x\neq 0$. Thus, 
\begin{equation*}
d=m2^{\left\lfloor \frac{n-1}{2}\right\rfloor }=m2^{\left( \frac{n-1}{2}%
\right) }
\end{equation*}
for some positive integer $m$, by Lemma \ref{dimensionLemma}. This proves
the statement since the dimension of $\mathbb{S}$ is $2^{\frac{n-1}{2}}$,
and the converse follows easily.
\end{proof}

The following result shows that nonsingular proper perturbations always
exist on Clifford modules over an odd dimensional manifold.

\begin{proposition}
\label{OddAlwaysExists}Suppose that the dimension $n$ of $M$ is odd. Let $E$
be a bundle of graded Clifford modules over $M$, and let $D$ be the
corresponding Dirac operator. Then there always exists a proper perturbation 
$Z$ of $D$; in particular the perturbation may be chosen to be invertible.
\end{proposition}

\begin{proof}
To prove that proper perturbations always exist, we simply take $\phi =i%
\mathbf{1}$ in Proposition \ref{SpinCOdd}.
\end{proof}

\begin{remark}
Note that the perturbation in the proof with $W^{\prime }=\mathbb{C}$ acts
on a bundle of rank $2^{\frac{n+1}{2}}$, so that the minimal rank in Theorem 
\ref{OddProperPerturbations} is sharp.
\end{remark}

\section{ Localization\label{LocalizationSection}}

In this section we will use \cite[Proposition 1.2]{Shu2} to study the
asymptotics of the spectrum of $D_{s}^{2}$ as $s\rightarrow \infty .$ (See
also \cite{Shu1}.)

Let 
\begin{eqnarray*}
H_{s} &:&=s^{-1}\left( D_{s}\right) ^{2}=s^{-1}D^{2}+ZD+DZ+sZ^{2} \\
&=&-s^{-1}A+B+sC,
\end{eqnarray*}
where $-A=D^{2}$ is a second order, elliptic, self-adjoint operator with a
non-negative principal symbol, and the operators $B=ZD+DZ$ and $C=Z^{2}$ are
self-adjoint bundle maps.

Each of the operators $H_{s},$ $A,$ $B,$ and $C$ has two self-adjoint
components, acting on $\Gamma \left( M,E^{+}\right) $ and on $\Gamma \left(
M,E^{-}\right) $, respectively. For example, 
\begin{equation*}
H_{s}^{+}:=s^{-1}\left. \left( D_{s}\right) ^{2}\right| _{\Gamma \left(
M,E^{+}\right) }=s^{-1}D_{s}^{-}D_{s}^{+}:\Gamma \left( M,E^{+}\right)
\rightarrow \Gamma \left( M,E^{+}\right) ,\text{ and }
\end{equation*}
\begin{equation*}
B^{-}:=\left. \left( ZD+DZ\right) \right| _{\Gamma \left( M,E^{-}\right)
}=Z^{+}D^{-}+D^{+}Z^{-}:\Gamma \left( M,E^{-}\right) \rightarrow \Gamma
\left( M,E^{-}\right) .
\end{equation*}

We now describe a model operator: a matrix harmonic oscillator that will
serve as an approximation of $\left( D_{s}\right) ^{2}$ near $\overline{x}$,
a singular point of $C$. We choose local coordinates $x=\left(
x_{1},...,x_{n}\right) $ and a trivialization of $E$ near $\overline{x}.$ We
assume that the volume associated to the Riemannian metric $g$ is the
Lebesgue volume element at the point $\overline{x}$ (this is easily done by
rescaling if necessary).

In the neighborhood of $\overline{x}$, operator $A$ becomes a $2m\times 2m$
block diagonal differential operator with two $m\times m$ blocks, where $m=%
\mathrm{rank}E^{+}=\mathrm{rank}E^{-}$. It has the form 
\begin{equation*}
A=\sum\limits_{1\leq i,j\leq n}A_{ij}\left( x\right) \frac{\partial ^{2}}{
\partial x_{i}\partial x_{j}}+A^{\left( 1\right) }.
\end{equation*}
Note that the operator $A^{\left( 1\right) }$ is at most first order. Let 
\begin{equation*}
A^{\left( 2\right) }=\sum\limits_{1\leq i,j\leq n}A_{ij}\left( \overline{x}
\right) \frac{\partial ^{2}}{\partial x_{i}\partial x_{j}};
\end{equation*}
i.e. $A^{\left( 2\right) }$ is the highest order part of $A$ taken at $x=%
\overline{x}$, it is a homogeneous second order differential operator with
constant coefficients that are $2m\times 2m$ block diagonal Hermitian
matrices.

We denote $\overline{B}=B\left( \overline{x}\right) ,$ so $\overline{B}$ is
just a $2m\times 2m$ block diagonal Hermitian matrix in the chosen
trivialization of $E$.

We also define 
\begin{equation*}
C^{\left( 2\right) }\left( x\right) =\frac{1}{2}\sum\limits_{1\leq i,j\leq n}%
\frac{\partial ^{2}C_{ij}}{\partial x_{i}\partial x_{j}}\left( \overline{x}
\right) x_{i}x_{j},
\end{equation*}
i.e. $C^{\left( 2\right) }$ is the quadratic part of the potential $C$ near $%
\overline{x}.$

We define the \textit{model operator} $K$ of $H_{s}$ at a singular point $%
\overline{x}$ to be 
\begin{equation*}
K\left( \overline{x}\right) =-A^{\left( 2\right) }+\overline{B}+C^{\left(
2\right) }.
\end{equation*}
We denote $m\times m$ blocks of $K\left( \overline{x}\right) $ as $K^{\pm
}\left( \overline{x}\right) .$ Each operator $K^{+}\left( \overline{x}%
\right) $ and $K^{-}\left( \overline{x}\right) $ has discrete spectrum,
since it is a quantum Hamiltonian of a $k$-dimensional harmonic oscillator 
\cite{Shu1}.

Let $\overline{x}_{1},\overline{x}_{2},...,\overline{x}_{N}$ be the list of
all singular points. Let 
\begin{equation*}
K^{\pm }=\bigoplus\limits_{l=1}^{N}K^{\pm }\left( \overline{x}_{l}\right)
\end{equation*}
be the model operators for $H_{s}^{\pm }$ on the set of all singular points
of $C=Z^{2}$. Denote the eigenvalues of $K^{\pm }$ by 
\begin{equation*}
\mu _{1}^{\pm }<\mu _{2}^{\pm }<\mu _{3}^{\pm }<...
\end{equation*}
and their multiplicities by 
\begin{equation*}
p_{1}^{\pm },p_{2}^{\pm },p_{3}^{\pm }....
\end{equation*}

\begin{proposition}
\label{ShubinTheorem}(Proposition 1.2 in \cite{Shu2}) If the $C\left(
x\right) \geq c\left\vert x-\overline{x}\right\vert ^{2}\mathbf{1}$, then
the eigenvalues of $H_{s}^{+}$ concentrate near the eigenvalues of the model
operator $K^{+}$. That is, for any positive integer $q$ there exists $%
s_{0}>0 $ and $c_{1}>0$ such that for any $s>$ $s_{0}$
\end{proposition}

\begin{theorem}
\begin{enumerate}
\item  there are precisely $p_{j}^{+}$ eigenvalues (multiplicities counted)
of $H_{s}^{+}$ in the interval 
\begin{equation*}
\left( \mu _{j}^{+}-c_{1}s^{-1/5},\mu _{j}^{+}+c_{1}s^{-1/5}\right)
,\,\,j=1,...,q;
\end{equation*}

\item  there no eigenvalues of $H_{s}^{+}$ in $\left( -\infty ,\mu
_{1}^{+}-c_{1}s^{-1/5}\right) $ and in the intervals 
\begin{equation*}
\left( \mu _{j}^{+}+c_{1}s^{-1/5},\mu _{j+1}^{+}-c_{1}s^{-1/5}\right)
,\,\,j=1,...,q;
\end{equation*}

\item  similar results also hold for operators $H_{s}^{-}$ and $K^{-}$.
\end{enumerate}
\end{theorem}

\begin{corollary}
\label{LocalizationCorollary}In the notation above, 
\begin{equation*}
\mathrm{ind}\left( D\right) =\dim \ker \left( K^{+}\right) -\dim \ker \left(
K^{-}\right) .
\end{equation*}
\end{corollary}

\begin{proof}
For each $s>0,$\ operators $H_{s}^{+}=s^{-1}$\ $D_{s}^{-}D_{s}^{+}$\ and $%
H_{s}^{-}=$\ $s^{-1}D_{s}^{+}D_{s}^{\_}$\ are positive elliptic self-adjoint
operators acting on sections of vector bundles over a compact smooth
manifold $M.$\ Therefore the operators $H_{s}^{+}$\ and $H_{s}^{-}$\ have
discrete spectra $\sigma \left( H_{s}^{\pm }\right) \subset \left[ 0,+\infty
\right) $ with finite multiplicities. By Proposition \ref{ShubinTheorem},
the spectra of $K^{+}$ and $K^{-}$ are also nonegative (and of course
discrete).

Choose any real number $r>0,$\ so that $r\ $is strictly less than the least
positive number in the union of the spectra of $K^{+}$\ and $K^{-}.$\ Then
for any $s>0$\ we have 
\begin{eqnarray*}
\mathrm{ind}\left( D\right) &=&\dim \ker \left( s^{-1}\left. \left(
D_{s}\right) ^{2}\right\vert _{\Gamma \left( M,E^{+}\right) }\right) -\dim
\ker \left( s^{-1}\left. \left( D_{s}\right) ^{2}\right\vert _{\Gamma \left(
M,E^{-}\right) }\right) , \\
&=&\dim \ker H_{s}^{+}-\dim \ker H_{s}^{-} \\
&=&\#\left\{ \sigma \left( H_{s}^{+}\right) \cap \left[ 0,r\right) \right\}
-\#\left\{ \sigma \left( H_{s}^{-}\right) \cap \left[ 0,r\right) \right\} ,
\end{eqnarray*}
because $D_{s}^{+}$ is an isomorphism between the eigenspaces of $H_{s}^{+}$
and of $H_{s}^{-}$ corresponding to nonzero eigenvalues. By choosing $s$\
sufficiently large in the formula above and applying Proposition \ref
{ShubinTheorem}, we obtain 
\begin{equation*}
\mathrm{ind}\left( D\right) =\dim \ker \left( K^{+}\right) -\dim \ker \left(
K^{-}\right) .
\end{equation*}
\end{proof}

\begin{remark}
\label{IndexZeroRemark}With the notation of Section \ref{preliminaryperturb}%
, if $Z:=\left( Z^{+},\left( Z^{+}\right) ^{\ast }\right) \in \Gamma \left(
M,\mathrm{End}\left( E^{+}\oplus E^{-}\right) \right) $ is everywhere
nonsingular and anticommutes with Clifford multiplication by vectors, then
the corollary implies that the index of the Dirac operator $D$ must be zero.
Proposition \ref{OddAlwaysExists} then yields a new proof that the index of
a Dirac operator (and thus any elliptic differential operator) on an
odd-dimensional manifold is zero.
\end{remark}

By applying Corollary \ref{LocalizationCorollary}, Proposition \ref{SpinC},
and the above remark to even dimensional manifolds, we obtain the following
corollary.

\begin{corollary}
\label{NoWIsomorphism}Let $M$ be even dimensional and spin$^{c}$. Let $%
E^{\pm }\cong \left( \mathbb{S}\otimes W\right) ^{\pm }=\left( \mathbb{S}%
^{+}\otimes W^{\pm }\right) \oplus \left( \mathbb{S}^{-}\otimes W^{\mp
}\right) $ be any graded, self-adjoint Clifford module over $M$. Then the
index of the Dirac operator corresponding to this Clifford module is zero if
the bundles $W^{+}$ and $W^{-}$ are isomorphic.
\end{corollary}

\begin{proof}
If there exists a bundle isomorphism $\phi ^{+}:W^{+}\rightarrow W^{-}$,
then $Z^{\pm }=\gamma \otimes \left( \phi ^{+},\left( \phi ^{+}\right)
^{*}\right) :E^{\pm }\rightarrow E^{\mp }$ is a bundle isomorphism that
anticommutes with Clifford multiplication. Therefore, the index is zero.
\end{proof}

\section{ Local Index Theory on $\mathbb{R}^{n}$\label{localIndexSection}}

In this section, if $P$ is a linear operator defined on a dense domain in a
Hilbert space, then $P^{*}$ denotes the formal adjoint of $P$.

Let $A_{1},...,A_{n},B_{1},...,B_{n}$ be $m\times m$ matrices, and let 
\begin{eqnarray*}
A &=&\sum_{j=1}^{n}A_{j}\partial _{j}\,,\,B\left( x\right)
=\sum_{j=1}^{n}x_{j}B_{j},\,\text{and} \\
Q &=&A+B\left( x\right) .
\end{eqnarray*}
We take the domain of $Q$ to be the space of compactly supported, smooth $%
\mathbb{C}^{m}$-valued functions on $\mathbb{R}^{n}$. We assume that

\begin{enumerate}
\item  $A$ is an elliptic operator.

\item  There is a positive constant $K$ such that $\left( B\left( x\right)
\right) ^{*}B\left( x\right) \geq K\left| x\right| ^{2}$for all $x\in %
\mathbb{R}^{n}$.

\item  For each $j$ and $k$, $A_{j}^{*}B_{k}-B_{k}^{*}A_{j}=0$.
\end{enumerate}

\begin{remark}
Note that the second condition above is equivalent to the fact that the
smallest eigenvalue of $\left( B\left( x\right) \right) ^{\ast }B\left(
x\right) $ is at least $K\left\vert x\right\vert ^{2}$. Thus, the same
inequality holds for $B\left( x\right) \left( B\left( x\right) \right)
^{\ast }$.
\end{remark}

Then 
\begin{eqnarray*}
Q^{\ast }Q &=&A^{\ast }A-\left( \sum_{j=1}^{n}A_{j}^{\ast }B_{j}\right)
+\left( B\left( x\right) \right) ^{\ast }B\left( x\right) \mathbf{,}\text{
and} \\
QQ^{\ast } &=&AA^{\ast }+\left( \sum_{j=1}^{n}A_{j}B_{j}^{\ast }\right)
+B\left( x\right) \left( B\left( x\right) \right) ^{\ast }
\end{eqnarray*}

By \cite[Theorem 2.2]{Ch}, the operator $\mathbf{Q}$ on $\mathbb{C}%
^{n}\times \mathbb{C}^{n}$ valued-functions defined by 
\begin{equation*}
\mathbf{Q}=\left( 
\begin{array}{cc}
0 & Q^{\ast } \\ 
Q & 0
\end{array}
\right) =\sum_{j=1}^{n}\left( 
\begin{array}{cc}
0 & -A_{j}^{\ast } \\ 
A_{j} & 0
\end{array}
\right) \partial _{j}+\left( 
\begin{array}{cc}
0 & \left( B\left( x\right) \right) ^{\ast } \\ 
B\left( x\right) & 0
\end{array}
\right)
\end{equation*}
is essentially self-adjoint. Its unique self-adjoint extension is defined as
the closure of \textbf{$Q$} with respect to the inner product norm $%
\left\Vert \cdot \right\Vert _{\mathbf{Q}}$ defined by $\left\Vert
u\right\Vert _{\mathbf{Q}}^{2}=\left\Vert u\right\Vert
_{L^{2}}^{2}+\left\Vert \mathbf{Q}u\right\Vert _{L^{2}}^{2}$.

Observe that 
\begin{equation*}
\mathbf{Q}^{2}=\left( 
\begin{array}{cc}
Q^{\ast }Q & 0 \\ 
0 & QQ^{\ast }
\end{array}
\right) .
\end{equation*}
Each of the operators $Q^{\ast }Q$ and $QQ^{\ast }$ are bounded from below
by 
\begin{equation*}
P=c\sum_{j=1}^{n}\ \left( -\partial _{j}^{2}+x_{j}^{2}\right) \mathbf{1}%
-\lambda \mathbf{1}
\end{equation*}
and above by 
\begin{equation*}
P^{\prime }=\frac{1}{c}\sum_{j=1}^{n}\ \left( -\partial
_{j}^{2}+x_{j}^{2}\right) \mathbf{1}+\lambda \mathbf{1}
\end{equation*}
for some constant $c>0$, and where $\lambda $ is the largest eigenvalue of $%
\left\vert \sum_{j=1}^{n}A_{j}B_{j}^{\ast }\right\vert $. Therefore, the
norm $\left\Vert \cdot \right\Vert _{\mathbf{Q}}$ defined above is
equivalent to the harmonic oscillator norm $\left\Vert \cdot \right\Vert $
defined by 
\begin{eqnarray*}
\left\Vert u\right\Vert ^{2} &=&\left\Vert u\right\Vert
_{L^{2}}^{2}+\sum_{j=1}^{n}\left\Vert \partial _{j}u\right\Vert
_{L^{2}}^{2}\ +\sum_{j=1}^{n}\left\Vert x_{j}u\right\Vert _{L^{2}}^{2} \\
&=&\left\Vert u\right\Vert _{L^{2}}^{2}+\sum_{j=1}^{n}\ \left\langle \left(
-\partial _{j}^{2}+x_{j}^{2}\right) u\left( x\right) ,u\left( x\right)
\right\rangle _{L^{2}}
\end{eqnarray*}
Denote by $\mathcal{H}$ the closure of the space of compactly supported,
smooth $\mathbb{C}^{m}$-valued functions on $\mathbb{R}^{n}$ with respect to
the norm $\left\Vert u\right\Vert $.

It is well known that $P+\tau \mathbf{1}$ has a compact inverse for some
constant $\tau $, so \textbf{$Q$}$^{2}+\tau \mathbf{1}$ must also have a
compact inverse. It follows that \textbf{$Q$}$^{2}$ has a finite dimensional
kernel. We conclude that both operators \textbf{$Q$} and $Q$ are Fredholm.

We summarize the arguments of this section in the following proposition.

\begin{proposition}
\label{localFredholmness}Suppose that $Q=\sum_{j=1}^{n}A_{j}\partial
_{j}\,+B\left( x\right) $ satisfies conditions (1) through (3) at the
beginning of this section. Then the closure of the elliptic operator $Q$ is
Fredholm on its domain $\mathcal{H}$.
\end{proposition}

\begin{corollary}
\label{matrixDeformationCorollary}Let $\left\{ \left. Q_{t}\,\right| t\in
\left[ 0,1\right] \,\right\} $ be a family of operators of the form 
\begin{equation*}
Q_{t}=\sum_{j=1}^{n}A_{j}\left( t\right) \partial
_{j}\,+\sum_{j=1}^{n}x_{j}B_{j}\left( t\right) ,
\end{equation*}
where $A_{j}\left( t\right) $ and $B_{j}\left( t\right) $ are continuous
families of matrices such that $\sum_{j=1}^{n}x_{j}A_{j}$ and $%
\sum_{j=1}^{n}x_{j}B_{j}$ are invertible for any nonzero $x\in \mathbb{R}%
^{n} $, satisfying 
\begin{equation*}
A_{j}^{*}B_{k}-B_{k}^{*}A_{j}=0
\end{equation*}
for each $j$ and $k$. Then the index of $Q_{t}$ is defined and is
independent of $t$.
\end{corollary}

\begin{proof}
It suffices to show that the family $\left\{ \left. Q_{t}:\mathcal{H}%
\rightarrow L^{2}\,\right\vert t\in \left[ 0,1\right] \,\right\} $ is a
continuous family in the norm topology. If $\max_{j}\left\vert A_{j}\left(
t_{1}\right) -A_{j}\left( t_{2}\right) \right\vert $ and $\max_{j}\left\vert
B_{j}\left( t_{1}\right) -B_{j}\left( t_{2}\right) \right\vert $ are both
less than $\delta >0$ (with respect to a fixed matrix norm), then for every $%
u\in \mathcal{H}$, 
\begin{eqnarray*}
&&\left\Vert \left( Q_{t_{1}}-Q_{t_{2}}\right) u\right\Vert _{L^{2}}^{2} \\
&=&\left\Vert \left( \sum_{j=1}^{n}\left( A_{j}\left( t_{1}\right)
-A_{j}\left( t_{2}\right) \right) \partial _{j}\,+\sum_{j=1}^{n}x_{j}\left(
B_{j}\left( t_{1}\right) -B_{j}\left( t_{2}\right) \right) \right)
u\right\Vert _{L^{2}}^{2} \\
&\leq &\delta ^{2}\left( \sum_{j=1}^{n}\left\Vert \partial _{j}u\right\Vert
_{L^{2}}^{2}\ +\sum_{j=1}^{n}\left\Vert x_{j}u\right\Vert
_{L^{2}}^{2}\right) \leq \delta ^{2}\left\Vert u\right\Vert ^{2}.
\end{eqnarray*}
\end{proof}

\section{ Local Calculations\label{LocalCalculationsSection}}

We now proceed with a calculation near the singular sets that evaluates the
index of the Dirac operator $D$. Suppose that $Z$ is a proper perturbation
of $D$ on $\Gamma \left( M,E\right) $ with a nonempty set of singular
points. This implies that there are restrictions on the rank of $E$ on the
associated grading and on the graded bundle map $Z$; see Section \ref
{globalClassificationSection} and Section \ref{properPerturbSection}. Let $U$
be the neighborhood of a singular point $\overline{x}$ as in Definition \ref
{properPerturbationDef}. We consider the operator $D_{s}=D+sZ$ on sections
of $E$, which has the local form 
\begin{equation*}
D_{s}=D+s\sum_{j}\left( x-\overline{x}\right) _{j}Z_{j},
\end{equation*}
where each $Z_{j}^{+}:E^{+}\rightarrow E^{-}$ is a locally defined,
self-adjoint bundle isomorphism with $Z_{j}=\left( Z_{j}^{+},\left(
Z_{j}^{+}\right) ^{*}\right) $.

We choose geodesic normal coordinates $x$ around $\overline{x}$ such that
the metric is the identity matrix at $\overline{x}$ and such that $\partial
_{j}$ commutes with $c\left( \partial _{k}\right) $ for all $j$ and $k$ at
the point $x=\overline{x}$. We wish to calculate the dimensions of $\ker
\left( \left. D_{s}^{2}\right| _{\Gamma \left( M,E^{\pm }\right) }\right) $.
By Corollary \ref{LocalizationCorollary}, it suffices to calculate the
dimensions of the kernels of the model operators $K^{\pm }\left( \overline{x}%
\right) $ corresponding to each singular point $\overline{x}$. In the
notation of the Section \ref{LocalizationSection}, the model operator is $%
K\left( \overline{x}\right) =-A^{\left( 2\right) }+\overline{B}+C^{\left(
2\right) }$, with 
\begin{eqnarray*}
A^{\left( 2\right) } &=&\sum_{j=1}^{n}\mathbf{1}\partial _{j}^{2}, \\
\overline{B} &=&\left. Z\circ D+D\circ Z\right| _{\overline{x}%
}=\sum_{j=1}^{n}c\left( \partial _{j}\right) Z_{j}\left( \overline{x}\right)
,\text{ and} \\
C^{\left( 2\right) } &=&\text{quadratic part of }Z^{2}=\left( \sum_{j}\left(
x-\overline{x}\right) _{j}Z_{j}\left( \overline{x}\right) \right) ^{2}
\end{eqnarray*}

Notice that 
\begin{equation*}
K\left( \overline{x}\right) =\left( 
\begin{array}{cc}
K^{+}\left( \overline{x}\right) & 0 \\ 
0 & K^{-}\left( \overline{x}\right)
\end{array}
\right) =\left( 
\begin{array}{cc}
D^{-}\left( \overline{x}\right) D^{+}\left( \overline{x}\right) & 0 \\ 
0 & D^{+}\left( \overline{x}\right) D^{-}\left( \overline{x}\right)
\end{array}
\right) ,
\end{equation*}
where the operator $D\left( \overline{x}\right) $ is defined by 
\begin{equation*}
D\left( \overline{x}\right) =\sum_{j}c\left( \partial _{j}\right) \partial
_{j}+\sum_{j}\left( x-\overline{x}\right) _{j}Z_{j}\left( \overline{x}%
\right) ,
\end{equation*}
which satisfies the hypothesis of Proposition \ref{localFredholmness}. Thus, 
$D\left( \overline{x}\right) $ is a Fredholm operator on $\mathbb{R}^{n}$,
and 
\begin{equation*}
\dim \ker \left( K^{+}\left( \overline{x}\right) \right) -\dim \ker \left(
K^{-}\left( \overline{x}\right) \right) =\mathrm{ind}_{\mathbb{R}^{n}}\left(
D\left( \overline{x}\right) \right) .
\end{equation*}

Corollary \ref{LocalizationCorollary} implies the following theorem.

\begin{theorem}
\label{masterpiece0}Let $D:\Gamma \left( M,E\right) \rightarrow \Gamma
\left( M,E\right) $ be the Dirac operator corresponding to a graded,
self-adjoint, Hermitian Clifford module over a closed manifold. Suppose that
there exists a proper perturbation (Definition \ref{properPerturbationDef}) $%
Z$ of $D$. Then near each singularity $\overline{x}$ of $Z$, we write $%
Z\left( x\right) =\sum_{j}\left( x-\overline{x}\right) _{j}Z_{j}$, where
each $Z_{j}^{+}:E^{+}\rightarrow E^{-}$ is a locally defined bundle
isomorphism with $Z_{j}=\left( Z_{j}^{+},\left( Z_{j}^{+}\right) ^{*}\right) 
$. Then the index of $D$ satisfies 
\begin{equation*}
\mathrm{ind}\left( D\right) =\sum_{\overline{x}\text{ singular}}\mathrm{ind}%
_{\mathbb{R}^{n}}\left( D\left( \overline{x}\right) \right) ,
\end{equation*}
where $D\left( \overline{x}\right) =\sum_{j}c\left( \partial _{j}\right)
\partial _{j}+\sum_{j}x_{j}Z_{j}\left( \overline{x}\right) $ is the operator
on $\mathbb{R}^{n}$ that maps $E_{\overline{x}}^{+}$-valued functions to $E_{%
\overline{x}}^{-}$-valued functions.
\end{theorem}

Let us now specialize to a case where we compute the local indices $\mathrm{%
ind}_{\mathbb{R}^{n}}\left( D\left( \overline{x}\right) \right) $ in terms
of the matrices $c\left( \partial _{1}\right) ,...,c\left( \partial
_{n}\right) ,Z_{1}\left( \overline{x}\right) ,...,Z_{n}\left( \overline{x}%
\right) $ directly. We will need to assume that $Z\left( x\right)
^{2}=q\left( x-\overline{x}\right) \mathbf{1}$, where $q$ is a positive
definite quadratic form and $\mathbf{1}$ is the identity map. The following
propositions show that there always exist local bundle maps with this
property. Furthermore, the coordinates may be chosen so that the $%
Z_{j}\left( \overline{x}\right) $ anticommute.

\begin{proposition}
Suppose the operator $Z\left( x\right) =\sum_{j=1}^{n}\left( x-\overline{x}%
\right) _{j}Z_{j}$ satisfies $Z\left( x\right) ^{2}=q\left( x-\overline{x}%
\right) \mathbf{1}$, where $q$ is a positive definite quadratic form and $%
\mathbf{1}$ is the identity map. Then there exist local coordinates $y$ and
Hermitian linear transformations $\widetilde{Z}_{2},...,\widetilde{Z}_{n}$
such that 
\begin{equation*}
Z\left( y\right) =\sum_{j=1}^{n}y_{j}\widetilde{Z}_{j}
\end{equation*}
and $\widetilde{Z}_{j}\widetilde{Z}_{k}+\widetilde{Z}_{k}\widetilde{Z}%
_{j}=2\delta _{jk}\mathbf{1}$. Furthermore, we have that $E\cong \mathbb{%
S\otimes }W$, where the rank of $W$ is a multiple of $2^{\left\lceil \frac{%
n-1}{2}\right\rceil }$.
\end{proposition}

\begin{proof}
If $Z\left( x\right) ^{2}=q\left( x-\overline{x}\right) \mathbf{1}$, then
there is a symmetric matrix $Q$, an orthogonal matrix $U,$ and a positive
diagonal matrix $D$ such that 
\begin{equation*}
q\left( x-\overline{x}\right) =Q\left( x-\overline{x}\right) \cdot \left( x-%
\overline{x}\right)
\end{equation*}
and such that $D=UQU^{T}$ is the identity. Let $y=\sqrt{D}U\left( x-%
\overline{x}\right) $. In the new coordinates we have 
\begin{equation*}
Z\left( y\right) =\sum_{j=1}^{n}y_{j}\widetilde{Z}_{j}
\end{equation*}
for the hermitian linear transformations $\widetilde{Z}_{j}=\sum_{k=1}^{n}%
\left( \sqrt{D}^{-1}U\right) _{jk}Z_{k}$. Then $\left( Z^{2}\right) \left(
y\right) =\left( \left\Vert y\right\Vert ^{2}\right) I$ implies the
following relations: 
\begin{equation*}
\widetilde{Z}_{j}\widetilde{Z}_{k}+\widetilde{Z}_{k}\widetilde{Z}%
_{j}=2\delta _{jk}.
\end{equation*}
Note then that $\left\{ i\widetilde{Z}_{j}\right\} $ becomes a set of
Clifford matrices, all of which commute with the given Clifford action. 
\textbf{\ }
\end{proof}

\begin{proposition}
Suppose that $r=\mathrm{rank}\left( E^{+}\right) $ is a multiple of $2^{n-1}$
. Then there exists a set of Hermitian, invertible linear transformations 
\newline
$\left\{ Z_{j}=\left( 
\begin{array}{ll}
0 & Z_{j}^{-} \\ 
Z_{j}^{+} & 0
\end{array}
\right) :=\left( 
\begin{array}{ll}
0 & \left( Z_{j}^{+}\right) ^{*} \\ 
Z_{j}^{+} & 0
\end{array}
\right) \right\} _{1\leq j\leq n}$ on the graded $\mathbb{C}\mathrm{l}\left( %
\mathbb{R}^{n}\right) $ module $E^{+}\oplus E^{-}\cong \mathbb{C}^{r}\oplus %
\mathbb{C}^{r}$ such that

\begin{enumerate}
\item  each $Z_{j}$ anticommutes with each $c\left( \partial _{k}\right) $
for $1\leq k\leq n$, and

\item  The operator $Z\left( x\right) =\sum_{j=1}^{n}\left( x-\overline{x}%
\right) _{j}Z_{j}$ satisfies $Z\left( x\right) ^{2}=q\left( x-\overline{x}
\right) \mathbf{1}$, where $q$ is a positive definite quadratic form and $%
\mathbf{1}$ is the identity map.
\end{enumerate}
\end{proposition}

\begin{proof}
Suppose that the rank of $E^{+}$ is a multiple of $2^{n-1}$, and $E$ is
endowed with a graded $\mathbb{C}\mathrm{l}\left( \mathbb{R}^{n}\right) $
action. Then there exists a graded $\mathbb{C}\mathrm{l}\left( \mathbb{R}%
^{2n}\right) $ action extending the original action on $E=E^{+}\oplus E^{-}$%
. If $\left\{ \beta _{j}\right\} _{1\leq j\leq n}$ is a set of generators
corresponding to Clifford multiplication on $E$ by the additional vectors,
then $\left\{ Z_{j}=i\beta _{j}\right\} _{1\leq j\leq n}$ is a set of
transformations that satisfy the conditions of the proposition.
\end{proof}

If the $\left\{ Z_{j}\left( \overline{x}\right) \right\} $ anticommute
(changing coordinates if necessary), then 
\begin{eqnarray*}
K\left( \overline{x}\right) &=&-\sum_{j=1}^{n}\partial
_{j}^{2}+\sum_{j=1}^{n}c\left( \partial _{j}\right) Z_{j}+\left(
\sum_{j}\left( x-\overline{x}\right) _{j}Z_{j}\right) ^{2} \\
&=&-\sum_{j=1}^{n}\partial _{j}^{2}+\sum_{j=1}^{n}c\left( \partial
_{j}\right) Z_{j}+\sum_{j=1}^{n}\left( x-\overline{x}\right)
_{j}^{2}Z_{j}^{2} \\
&=&\sum_{j=1}^{n}\left( -\partial _{j}^{2}+c\left( \partial _{j}\right)
Z_{j}+\left( x-\overline{x}\right) _{j}^{2}\left( c\left( \partial
_{j}\right) Z_{j}\right) ^{2}\right) ,
\end{eqnarray*}
where each $Z_{j}$ is evaluated at $x=\overline{x}$.

Observe that the operators $L_{j}=c\left( \partial _{j}\right) Z_{j}$ are
Hermitian and commute with each other and thus can be diagonalized
simultaneously. Let $v$ be a common eigenvector of each operator $L_{j}$
corresponding to the eigenvalue $\lambda _{j}$. Letting $f$ be a function of 
$x$, we have 
\begin{equation*}
K\left( \overline{x}\right) \left( fv\right) =\left( \sum_{j=1}^{n}\left(
-\partial _{j}^{2}+\lambda _{j}+\lambda _{j}^{2}\left( x-\overline{x}\right)
_{j}^{2}\right) f\right) v.
\end{equation*}
The section $fv$ is in the kernel of $K\left( \overline{x}\right) $ if and
only if each $\lambda _{j}$ is negative, and, up to a constant, $fv=\exp
\left( \frac{1}{2}\sum_{j}\lambda _{j}\left( x-\overline{x}\right)
_{j}^{2}\right) v$. Thus the dimension of the kernel of $K\left( \overline{x}%
\right) $ is the dimension of the intersection of the direct sum of
eigenspaces $E_{\lambda }\left( L_{j}\right) $ of $L_{j}=c\left( \partial
_{j}\right) Z_{j}$ corresponding to negative eigenvalues. Note that $L_{j}$
maps $E^{+}$ to itself (call the restriction $L_{j}^{+}$), so that the
dimension of the kernel of $D_{s}^{+}$ restricted to this neighborhood is
simply the dimension of $\bigcap_{j}\left( \bigoplus_{\lambda <0}E_{\lambda
}\left( L_{j}^{+}\right) \right) $.

The calculation above and Corollary \ref{LocalizationCorollary} imply the
following theorem.

\begin{theorem}
\label{Masterpiece}Let $D:\Gamma \left( M,E\right) \rightarrow \Gamma \left(
M,E\right) $ be the Dirac operator corresponding to a graded, self-adjoint,
Hermitian Clifford module over a closed manifold. Suppose that there exists
a proper perturbation (Definition \ref{properPerturbationDef}) $Z$ of $D$.
Then near each singularity $\overline{x}$ of $Z$, we write $Z\left( x\right)
=\sum_{j}\left( x-\overline{x}\right) _{j}Z_{j}$, where each $%
Z_{j}^{+}:E^{+}\rightarrow E^{-}$ is a locally defined bundle isomorphism
with $Z_{j}=\left( Z_{j}^{+},\left( Z_{j}^{+}\right) ^{*}\right) $. We
assume that we may choose $Z$ such that $Z_{j}\left( \overline{x}\right)
Z_{k}\left( \overline{x}\right) =-Z_{k}\left( \overline{x}\right)
Z_{j}\left( \overline{x}\right) $ for all $j\neq k$. Define the Hermitian
linear transformations 
\begin{equation*}
L_{j}^{\pm }\left( \overline{x}\right) =\left. c\left( \partial _{j}\right)
Z_{j}\right| _{E_{\overline{x}}^{\pm }}.
\end{equation*}
Then 
\begin{equation*}
\mathrm{ind}\left( D\right) =\sum_{\overline{x}}\left( \dim \left[
\bigcap_{j}\left( \bigoplus_{\lambda <0}E_{\lambda }\left( L_{j}^{+}\left( 
\overline{x}\right) \right) \right) \right] -\dim \left[ \bigcap_{j}\left(
\bigoplus_{\lambda <0}E_{\lambda }\left( L_{j}^{-}\left( \overline{x}\right)
\right) \right) \right] \right) ,
\end{equation*}
where the sum is taken over all the singular points $\overline{x}$ of $Z$.
\end{theorem}

\begin{remark}
The assumption that $Z_{j}\left( \overline{x}\right) Z_{k}\left( \overline{x}
\right) =-Z_{k}\left( \overline{x}\right) Z_{j}\left( \overline{x}\right) $
for all $j\neq k$ is natural, since it appears in all local calculations
where researchers have used Witten deformation (see papers mentioned in the
introduction). In fact, a stronger condition, that $Z_{j}^{2}$ is a scalar
multiplication by a function and that the anticommutivity holds for all $x$
near $\overline{x}$, appears in all previous work.
\end{remark}

We now apply Theorem \ref{Masterpiece} when $n=\dim M$ is even. Suppose that
the hypotheses of Theorem \ref{EvenProperPerturbations} hold, and suppose
that $Z$ has a nonempty set of singular points. Let $U$ be the neighborhood
of a singular point $\overline{x}$ as in Definition \ref
{properPerturbationDef}. Then the operator $D_{s}=D+sZ$ over $E|_{U}\cong %
\mathbb{S|}_{U}\otimes W$ has the form 
\begin{eqnarray*}
D_{s} &=&D+s\gamma \otimes \phi \\
&=&D+s\gamma \otimes \sum_{j}\left( x-\overline{x}\right) _{j}\phi _{j},
\end{eqnarray*}
with notation as in Definition \ref{properPerturbationDef} and Proposition 
\ref{SpinC}, where each $\phi _{j}^{+}:W^{+}\rightarrow W^{-}$ is a locally
defined self-adjoint bundle isomorphism with $\phi _{j}=\left( \phi
_{j}^{+},\left( \phi _{j}^{+}\right) ^{*}\right) $. The theorem becomes:

\begin{corollary}
\label{EvenMasterpiece}Let $D:\Gamma \left( M,E\right) \rightarrow \Gamma
\left( M,E\right) $ be the Dirac operator corresponding to a graded,
self-adjoint, Hermitian Clifford module over an even-dimensional manifold.
Let $Z$ be a proper perturbation of $D$. Then near each singularity $%
\overline{x}$ of $Z$, we write $E^{\pm }\cong \left( S^{\pm }\otimes W^{\pm
}\right) \oplus \left( S^{\pm }\otimes W^{\pm }\right) $ and $Z\left(
x\right) =\sum_{j}\left( x-\overline{x}\right) _{j}\gamma \otimes \phi _{j}$%
, where each $\phi _{j}^{+}:W^{+}\rightarrow W^{-}$ is a locally defined
bundle isomorphism with $\phi _{j}=\left( \phi _{j}^{+},\left( \phi
_{j}^{+}\right) ^{*}\right) $. Assume that $\phi _{j}\phi _{k}=-\phi
_{k}\phi _{j}$ for all $j\neq k$. Define the Hermitian linear
transformations 
\begin{equation*}
L_{j}^{\pm }\left( \overline{x}\right) =\left. \left( \left( c\left(
\partial _{j}\right) \gamma \right) \otimes \phi _{j}\right) \right| _{E_{%
\overline{x}}^{\pm }}.
\end{equation*}
Then 
\begin{equation}
\mathrm{ind}\left( D\right) =\sum_{\overline{x}}\left( \dim \left[
\bigcap_{j}\left( \bigoplus_{\lambda <0}E_{\lambda }\left( L_{j}^{+}\left( 
\overline{x}\right) \right) \right) \right] -\dim \left[ \bigcap_{j}\left(
\bigoplus_{\lambda <0}E_{\lambda }\left( L_{j}^{-}\left( \overline{x}\right)
\right) \right) \right] \right) ,  \label{Evenformula}
\end{equation}
where the sum is taken over all the singular points $\overline{x}$ of $Z$.
\end{corollary}

\begin{remark}
\label{OddMasterpiece}The proof is easily modified for the odd-dimensional
case. Locally we write $E|_{U}\cong \left. S\otimes \left( W^{\prime }\oplus
W^{\prime }\right) \right| _{U}$ as in Proposition \ref{SpinCOdd}. Then $%
Z=\sum_{j}\left( x-\overline{x}\right) _{j}\mathbf{1}\otimes \left( 
\begin{array}{ll}
0 & \phi _{j} \\ 
-\phi _{j} & 0
\end{array}
\right) $, where each $\phi _{j}:W^{\prime }\rightarrow W^{\prime }$ is a
locally defined, skew-adjoint bundle isomorphism. Formula (\ref{Evenformula}
) is valid with the new $L_{j}^{\pm }\left( \overline{x}\right) $ given by 
\begin{equation*}
L_{j}^{\pm }\left( \overline{x}\right) =\left. c^{+}\left( \partial
_{j}\right) \otimes \left( 
\begin{array}{ll}
0 & \phi _{j} \\ 
\phi _{j} & 0
\end{array}
\right) \right| _{E_{\overline{x}}^{\pm }}.
\end{equation*}
In this case, the formula shows that the sum of the local indices is always
zero, since $\mathrm{ind}\left( D\right) =0$ in odd dimensions.
\end{remark}

\section{Examples\label{ExamplesSection}}

\subsection{ The de Rham operator\label{deRham}}

The bundle $E=\Lambda ^{\bullet }T^{*}M\otimes \mathbb{C}$ of complex-valued
forms is a left Clifford module with the canonical Clifford action defined
by 
\begin{equation*}
l(v)\omega =v^{\flat }\wedge \omega -v\lrcorner \omega ,\;v\in
T_{x}M,\;\omega \in \Lambda ^{\bullet }T_{x}^{*}M\otimes \mathbb{C}.
\end{equation*}
Here $v^{\flat }$ denotes the covector dual to $v$ and $v\lrcorner $ is the
contraction with vector $v$. Similarly, $E$ is also a right Clifford module
with the canonical right Clifford action on $\Lambda ^{p}T^{*}M\otimes %
\mathbb{C}$ defined by 
\begin{equation*}
r(v)\omega =(-1)^{p}\left( v^{\flat }\wedge \omega +v\lrcorner \omega
\right) .
\end{equation*}
The corresponding Dirac operator is the \textit{de Rham operator} 
\begin{equation*}
D=\sum l\left( e_{j}\right) \nabla _{e_{j}}=d+d^{*},
\end{equation*}
where $\left\{ e_{1},...,e_{n}\right\} $ is a local orthonormal basis of $TM$%
.

Let the dimension of $M$ be even. Let $\mathbb{S}$ denote a $\mathrm{spin}%
^{c}$ bundle, which always exists locally but may not exist globally if $M$
is not a $\mathrm{spin}^{c}$ manifold. Let $\mathbb{S}^{\ast }$ denote the
dual bundle to $\mathbb{S}$. Using a bundle metric on $\mathbb{S}$, identify 
$\mathbb{S}^{\ast }$ and $\overline{\mathbb{S}}.$ The result below follows
from the representation theory of Clifford algebras.

\begin{lemma}
\label{Isomorphism}(follows from Ch. IV in \cite{L-M}) Through the
isomorphism of Clifford modules $\mathbb{S}\otimes \mathbb{S}^{*}\cong
\Lambda ^{\bullet }T^{*}M\otimes \mathbb{C}$, the corresponding Clifford
actions by vectors are intertwined in the following commutative diagrams: 
\begin{equation*}
\begin{array}{lll}
\mathbb{S}\otimes \mathbb{S}^{*} & \overset{\cong }{\longrightarrow } & 
\Lambda ^{\bullet }T^{*}M\otimes \mathbb{C} \\ 
\downarrow {c}\left( v\right) {\otimes 1} & \circlearrowleft  & \downarrow {l%
}\left( v\right)  \\ 
\mathbb{S}\otimes \mathbb{S}^{*} & \overset{\cong }{\longrightarrow } & 
\Lambda ^{\bullet }T^{*}M\otimes \mathbb{C}
\end{array}
\,\,\,\,\,\,\,\,\,\,\,\,\,\,\,\,\,
\begin{array}{lll}
\mathbb{S}\otimes \mathbb{S}^{*} & \overset{\cong }{\longrightarrow } & 
\Lambda ^{\bullet }T^{*}M\otimes \mathbb{C} \\ 
\downarrow {1\otimes c}\left( -v\right)  & \circlearrowleft  & \downarrow {r}%
\left( v\right)  \\ 
\mathbb{S}\otimes \mathbb{S}^{*} & \overset{\cong }{\longrightarrow } & 
\Lambda ^{\bullet }T^{*}M\otimes \mathbb{C}
\end{array}
\end{equation*}
\end{lemma}

The grading of $E$ into even and odd forms $E=E^{+}\oplus E^{-}$, where $%
E^{+}=\Lambda ^{\text{\textrm{even}}}T^{\ast }M\otimes \mathbb{C}$ and $%
E^{-}=\Lambda ^{\text{\textrm{odd}}}T^{\ast }M\otimes \mathbb{C}$, is not
natural; it does not come from the action of the chirality operator $\gamma $
on $E.$ Under the isomorphism \textbf{\ }$E=\Lambda ^{\bullet }T^{\ast
}M\otimes \mathbb{C}\cong \mathbb{S}\otimes \mathbb{S}^{\ast }$ we have $%
E^{+}=\left( \mathbb{S}^{+}\otimes \left( \mathbb{S}^{+}\right) ^{\ast
}\right) \oplus \left( \mathbb{S}^{-}\otimes \left( \mathbb{S}^{-}\right)
^{\ast }\right) $ and $E^{-}=\left( \mathbb{S}^{+}\otimes \left( \mathbb{S}%
^{-}\right) ^{\ast }\right) \oplus \left( \mathbb{S}^{-}\otimes \left( %
\mathbb{S}^{+}\right) ^{\ast }\right) $. Since the Dirac operator only acts
on the first component of $\mathbb{S}\otimes \mathbb{S}^{\ast }$, we can
deform $D$ by any linear operator $Z=\gamma \otimes \phi ,$ where $\gamma $
is the chirality operator and $\phi :\mathbb{S}^{\pm }\rightarrow \mathbb{S}%
^{\mp }$ is a bundle map as in Proposition \ref{SpinC}.

The following example of such a deformation is useful in the proof of the
Poincar\'{e}-Hopf Theorem (\ref{PoincareHopf}). Let $V$ be a smooth vector
field on $M$, then for each $p$-form $\omega $ define 
\begin{equation*}
Z_{V}\omega =\left( V^{\flat }\wedge +V\lrcorner \right) \omega =\left(
-1\right) ^{p}r\left( V\right) \omega :E^{\pm }\rightarrow E^{\mp }.
\end{equation*}

Recall that a point $\overline{x}\in M$ is called \textit{singular (or
critical )} point of $V\left( x\right) =\sum_{k}V_{k}\left( x\right)
\partial _{k}$ if for all $k=1,...,n,$ $V_{k}(\overline{x})=0$. A singular
point of $V$ is called \textit{non-degenerate} if $\det \left( \partial
V_{k}/\partial x_{i}\right) (\overline{x})\neq 0.$ The index of a singular
point $\overline{x}$ is defined as follows 
\begin{equation*}
\text{ind}_{V}\left( \overline{x}\right) =\text{sign }\det \left( \partial
V_{k}/\partial x_{i}\right) (\overline{x}).
\end{equation*}

\begin{remark}
Note that $\overline{x}\in M$ is a nondegenerate singular point of $V$ if
and only if $\overline{x}$ is a proper singular point of the endomorphism $%
Z_{V}$ on forms.
\end{remark}

Observe that the map $\omega \longmapsto \left( -1\right) ^{p}\omega $ on $p$
-forms $\omega $ is given by the map $\gamma \otimes \gamma $ on $\mathbb{S}%
\otimes \mathbb{S}^{\ast }=\mathbb{S}\otimes \overline{\mathbb{S}}$, since $%
\gamma $ is the identity (respectively, minus the identity) on even
(respectively, odd) spinors. Thus, 
\begin{equation*}
Z_{V}=\left( -1\right) ^{p}r\left( V\right) =\left( \mathbf{1}\otimes
c\left( V\right) \right) \circ \left( \gamma \otimes \gamma \right) =\left(
\gamma \otimes c\left( V\right) \gamma \right) .
\end{equation*}
so that $Z_{V}=\gamma $ $\otimes \phi $ with $\phi =c\left( V\right) \gamma
\ $.\thinspace

Near each nondegenerate singular point $\overline{x}$ of the vector field $V$%
, we choose local coordinates centered at $\overline{x}$ so that $V\left(
x\right) =\sum_{j=1}^{n}x_{j}\sum_{k=1}^{n}V_{jk}\partial _{k}+\mathcal{O}%
\left( \left| x\right| ^{2}\right) $, where $\left( V_{jk}\right) $ is an
invertible matrix. The model operator $K\left( \overline{x}\right) $ from
Section \ref{LocalCalculationsSection} is 
\begin{equation*}
K\left( \overline{x}\right) =-\sum_{j=1}^{n}\mathbf{1}\partial
_{j}^{2}+\sum_{j,k=1}^{n}V_{jk}\left( c\left( \partial _{j}\right) \gamma
\otimes c\left( \partial _{k}\right) \gamma \right)
+\sum_{j,k,m=1}^{n}x_{j}x_{k}\left( VV^{T}\right) _{jk}.
\end{equation*}
We observe that $K\left( \overline{x}\right) =D\left( \overline{x}\right)
^{2}$ on even forms, where 
\begin{equation*}
D\left( \overline{x}\right) :\Gamma \left( \mathbb{R}^{n},\Lambda ^{\text{%
\textrm{even/odd}}}T^{*}\mathbb{R}^{n}\otimes \mathbb{C}\right) \rightarrow
\Gamma \left( \mathbb{R}^{n},\Lambda ^{\text{\textrm{odd/even}}}T^{*}%
\mathbb{R}^{n}\otimes \mathbb{C}\right)
\end{equation*}
is defined by 
\begin{equation*}
D\left( \overline{x}\right) =\sum_{j=1}^{n}\left( c\left( \partial
_{j}\right) \otimes \mathbf{1}\right) \partial
_{j}+\sum_{j,k=1}^{n}x_{j}V_{jk}\left( \gamma \otimes c\left( \partial
_{k}\right) \gamma \right) .
\end{equation*}

By Theorem \ref{masterpiece0}, 
\begin{equation}
\chi \left( M\right) =\mathrm{ind}\left( d+d^{*}\right) =\sum_{V\left( 
\overline{x}\right) =0}\mathrm{ind}_{\mathbb{R}^{n}}\left( D\left( \overline{%
x}\right) \right) .  \label{indexEquationLocalVF}
\end{equation}
By Corollary \ref{matrixDeformationCorollary}, each of the integers $\mathrm{%
ind}_{\mathbb{R}^{n}}\left( D\left( \overline{x}\right) \right) $ does not
change if the vector field $V$ is deformed continuously while remaining
nondegenerate and without introducing additional zeros. Every vector field
on $\mathbb{R}^{n}$ with a single, nondegenerate zero at the origin may be
continously deformed in this way to $V\left( x\right) =\pm x_{1}\partial
_{1}+\sum_{j=2}^{n}x_{j}\partial _{j}$, depending only on the index $\pm 1$.
Thus, to evaluate the right hand side of Equation (\ref{indexEquationLocalVF}%
), it suffices to calculate the index of the corresponding operator 
\begin{equation*}
D\left( \overline{x}\right) =\sum_{j=1}^{n}\left( c\left( \partial
_{j}\right) \otimes \mathbf{1}\right) \partial _{j}\pm x_{1}\left( \gamma
\otimes c\left( \partial _{1}\right) \gamma \right)
+\sum_{j=2}^{n}x_{j}\left( \gamma \otimes c\left( \partial _{j}\right)
\gamma \right) .
\end{equation*}

The operator above satisfies the conditions of Theorem \ref{Masterpiece}.
Note that the local zeroth order operator $Z$ corresponds to 
\begin{eqnarray*}
Z_{V} &=&\gamma \otimes \phi =\gamma \otimes \left[ \pm x_{1}c\left(
\partial _{1}\right) \gamma +\sum_{j=2}^{n}x_{j}\left( c(\partial
_{j})\gamma \right) \right] \\
&=&\gamma \otimes \left[ \sum_{j=1}^{n}x_{j}\phi _{j}\right] .
\end{eqnarray*}
In the notation of Section \ref{LocalCalculationsSection} we have 
\begin{eqnarray*}
L_{1} &:&=\pm \left( c\left( \partial _{1}\right) \gamma \right) \otimes
\left( c\left( \partial _{1}\right) \gamma \right) \\
L_{j} &:&=\left( c\left( \partial _{j}\right) \gamma \right) \otimes \left(
c\left( \partial _{j}\right) \gamma \right) \,\,\text{for }j\geq 2
\end{eqnarray*}
If $\omega $ is a $p$-form, then for instance 
\begin{eqnarray*}
L_{j}\omega &=&(-1)^{p}l\left( \partial _{j}\right) r\left( \partial
_{j}\right) \omega \\
&=&\left( dx_{j}\wedge -dx_{j}\lrcorner \right) \left( dx_{j}\wedge
+dx_{j}\lrcorner \right) \omega \text{ for }j\geq 2.
\end{eqnarray*}
If, in addition, $\omega =dx_{i_{1}}\wedge ...\wedge dx_{i_{p}},$ then 
\begin{equation*}
L_{j}\omega =\left\{ 
\begin{array}{l}
\omega ,\text{ if }j\in \left\{ i_{1},...,i_{p}\right\} \ \text{and }j\geq 2
\\ 
-\omega ,\text{ if }j\notin \left\{ i_{1},...,i_{p}\right\} \ \text{and }%
j\geq 2 \\ 
\pm \omega ,\text{ if }j=1\in \left\{ i_{1},...,i_{p}\right\} \\ 
\mp \omega ,\text{ if }j=1\notin \left\{ i_{1},...,i_{p}\right\}
\end{array}
\right.
\end{equation*}
The proposition below follows.

\begin{proposition}
\label{NegativeEigenspace}In a neighborhood of a zero of index $\pm 1$, then
the form $\omega $ belongs to the negative eigenspace of $L_{j}$ for each $j$
if and only if 
\begin{equation*}
\omega =\left\{ 
\begin{array}{l}
c\text{ if }\mathrm{ind}_{V}\left( \overline{x}\right) =1 \\ 
cdx_{1}\text{ if }\mathrm{ind}_{V}\left( \overline{x}\right) =-1
\end{array}
\right. .
\end{equation*}
for some constant $c$.
\end{proposition}

Since $\left( \Lambda ^{\bullet }T^{*}\mathbb{R}^{n}\otimes \mathbb{C}
\right) ^{+}=\left( \Lambda ^{\mathrm{even}}T^{*}\mathbb{R}^{n}\otimes %
\mathbb{C}\right) $, the local indices of $D^{\pm }\left( \overline{x}%
\right) $ --- and thus the index $\mathrm{ind}_{\mathbb{R}^{n}}\left(
D\left( \overline{x}\right) \right) $ --- are now clearly determined.

\begin{corollary}
We have $\mathrm{ind}_{\mathbb{R}^{n}}\left( D\left( \overline{x}\right)
\right) =\mathrm{ind}_{V}\left( \overline{x}\right) $.
\end{corollary}

The Poincar\'{e}-Hopf Theorem now follows from Equation (\ref
{indexEquationLocalVF}).

\begin{theorem}
(Poincar\'{e}-Hopf Theorem) Let $V$ be a smooth vector field with
non-degenerate singular points on an even-dimensional, smooth manifold $M$.
Let $n^{\pm }$ denote the number of singular points of $V$ with the index $%
\pm 1.$ Then 
\begin{equation*}
\chi \left( M\right) =n^{+}-n^{-},
\end{equation*}
where $\chi \left( M\right) $ is the Euler characteristic of $M.$
\end{theorem}

\begin{remark}
The Poincar\'{e}-Hopf Theorem on odd-dimensional manifolds can be proved in
a similar way.
\end{remark}

\begin{corollary}
\label{NoSpinIso}The Euler characteristic of a closed, even-dimensional spin$%
^{c}$ manifold is zero if and only if the spin$^{c}$ bundles $\mathbb{S}^{+}$
and $\mathbb{S}^{-}$ are isomorphic. The statement is true for any choice of
spin$^{c}$ bundle $\mathbb{S}$.
\end{corollary}

\begin{proof}
If the Euler characteristic of the manifold is zero, then there exists a
nonzero vector field. Clifford multiplication by this vector field provides
the needed isomorphism. On the other hand, if a bundle isomorphism $\psi :%
\mathbb{S}^{+}\rightarrow \mathbb{S}^{-}$ does exist, it induces a map $\phi
:\mathbb{S}^{\ast }\rightarrow \mathbb{S}^{\ast }$ defined using $\phi
^{+}=\psi ^{\ast }:\left( \mathbb{S}^{+}\right) ^{\ast }\rightarrow \left( %
\mathbb{S}^{-}\right) ^{\ast }$ and $\phi ^{-}=\left( \phi ^{+}\right)
^{\ast }$ and thus a proper perturbation $Z=\gamma \otimes \phi $ of the de
Rham operator with no singular points.
\end{proof}

\subsection{ The Euler characteristic and sections of the conformal Pin
bundle}

In this section, we again assume that the dimension $n$ of $M$ is even. We
will show that the Euler characteristic is the sum of indices of a
nondegenerate section of the conformal Pin bundle.

Consider the subset $P_{x}M$ of all elements of the $\mathbb{C}\mathrm{l}%
\left( T_{x}M\right) $ of the form 
\begin{equation*}
v_{1}v_{2}...v_{r}
\end{equation*}
where $v_{j}\in T_{x}M$ for each $j$. This set forms a monoid, and it is the
same as 
\begin{equation*}
P_{x}M=\left\{ \lambda \alpha \,|\,\alpha \in \mathrm{Pin}\left(
T_{x}M\right) ,\lambda \in \mathbb{R}\right\} \subset \mathbb{C}\mathrm{l}
\left( T_{x}M\right)
\end{equation*}
(see \cite[p. 12ff]{L-M}). Note that $P_{x}M-\left\{ 0\right\} $ is the
conformal Pin group of $T_{x}M$. Let $PM$ denote the corresponding bundle of
monoids over $M$. Let $P^{+}M$ and $P^{-}M$ be defined by 
\begin{equation*}
P^{\pm }M=PM\cap \mathbb{C}\mathrm{l}^{\pm }\left( T_{x}M\right) .
\end{equation*}

We say that $\overline{x}\in M$ is a \emph{\ nondegenerate zero} of $\beta
\in \Gamma \left( P^{-}M\right) $ if on a sufficiently small neighborhood $U$
of $\overline{x}$ we have

\begin{enumerate}
\item  $\beta |_{\overline{x}}=0$, and

\item  in local coordinates $x$ on $U$ , there exist invertible $\beta
_{j}\in \Gamma \left( \left. P^{-}M\right| _{U}\right) $ for $1\leq j\leq %
n=\dim M$ over $U$ such that $\beta =\sum_{j}\left( x-\overline{x}\right)
_{j}\beta _{j}$ on $U$, and $\beta $ is invertible on $U\setminus \left\{ 
\overline{x}\right\} $.
\end{enumerate}

If $\overline{x}\in M$ is a nondegenerate zero of $\beta \in \Gamma \left(
P^{-}M\right) $, then on some neighborhood $U$ of $\overline{x}$, 
\begin{equation*}
\left. \beta \right\vert _{U}=\left. W_{1}W_{2}...W_{r}\right\vert _{U}
\end{equation*}
for some vector fields $W_{1},...,W_{r}$. Since $\beta $ is nondegenerate,
only one of the vector fields (say $W_{j}$) is zero at $\overline{x}$, and $%
\overline{x}$ is a nondegenerate zero of $W_{j}$.

\begin{lemma}
\label{CommutingLemma}For any two vector fields $B_{1}$ and $B_{2}$ such
that $B_{2}$ is nonzero, 
\begin{equation*}
B_{2}B_{1}=\widetilde{B_{1}}B_{2},
\end{equation*}
where $\widetilde{B_{1}}$ is the vector field defined by 
\begin{equation*}
\widetilde{B_{1}}=B_{1}^{\Vert }-B_{1}^{\bot },
\end{equation*}
where $B_{1}^{\Vert }$ and $B_{1}^{\bot }$ are the components of $B_{1}$ in
directions parallel to and perpendicular to $B_{2}$. If $B_{1}$ has a
nondegenerate zero at a point where $B_{2}$ is nonzero, the index of $%
\widetilde{B_{1}}$ at the point is the opposite of the index of $B_{1}$ at
the point.
\end{lemma}

\begin{proof}
The equation follows from the construction of $\widetilde{B_{1}}$. The
vector field $-\widetilde{B_{1}}$ is the reflection of $B_{1}$ in the plane
perpendicular to $B_{2}$. Since the dimension is even, the index of $%
\widetilde{B_{1}}$ is the opposite of the index of $B_{1}$.
\end{proof}

Lemma \ref{CommutingLemma} implies that $\overline{x}\in M$ is a
nondegenerate zero of $\beta \in \Gamma \left( P^{-}M\right) $ if and only
if in local coordinates $x$ on a sufficiently small neighborhood $U$ of $%
\overline{x}$ , there exists a locally defined vector field $V_{1}$ with an
isolated nondegenerate zero at $\overline{x}$ and a collection of nonzero
vector fields $V_{2},...,V_{r}$ such that 
\begin{equation*}
\left. \beta \right\vert _{U}=\left. V_{1}V_{2}...V_{r}\right\vert _{U}.
\end{equation*}

\begin{definition}
Given a section $\beta \in \Gamma \left( P^{-}M\right) $ with a
nondegenerate zero at $\overline{x}\in M$, we define the \emph{index }$%
\mathrm{ind}_{\beta }\left( \overline{x}\right) $of\emph{\ }$\beta $ at $%
\overline{x}$ to be the index $\mathrm{ind}_{V_{1}}\left( \overline{x}%
\right) $of any vector field $V_{1}$ at $\overline{x}$ such that in a
neighborhood $U$ of $\overline{x}$, 
\begin{equation*}
\left. \beta \right| _{U}=\left. V_{1}V_{2}...V_{r}\right| _{U},
\end{equation*}
where $V_{1}$ has an isolated nondegenerate zero at $\overline{x}$ and the
vector fields $V_{2},...,V_{r}$ are nonzero at $\overline{x}$.
\end{definition}

\begin{lemma}
Given a section $\beta \in \Gamma \left( P^{-}M\right) $ with a
nondegenerate zero at $\overline{x}\in M$, the index\emph{\ }of\emph{\ }$
\beta $ at $\overline{x}$ is well-defined.
\end{lemma}

\begin{proof}
Suppose that we are given two different local expressions for $\beta $ on a
sufficiently small neighborhood $U$ of $\overline{x}$: 
\begin{eqnarray*}
\left. \beta \right| _{U} &=&\left. V_{1}V_{2}...V_{r}\right| _{U} \\
&=&\left. W_{1}W_{2}...W_{r^{\prime }}\right| _{U},
\end{eqnarray*}
where $V_{1}$ and $W_{1}$ have isolated nondegenerate zeros at $\overline{x}%
=0$ and the vector fields $V_{2},...,V_{r},W_{2},...,W_{r^{\prime }}$ are
nonzero at $0$. Without loss of generality, we replace $V_{1}$ and $W_{1}$
with their linear parts 
\begin{equation*}
V_{1}=\sum_{j=1}^{n}x_{j}V_{1j},W_{1}=\sum_{j=1}^{n}x_{j}W_{1j},
\end{equation*}
and we replace $V_{2},...,V_{r},W_{2},...,W_{r^{\prime }}$ with their values
at $0$. Thus, the equation for $\beta $ above implies that 
\begin{equation*}
V_{1j}V_{2}...V_{r}=W_{1j}W_{2}...W_{r^{\prime }}.
\end{equation*}
We multiply on the right by the inverse of $V_{2}...V_{r}$ to obtain 
\begin{eqnarray*}
V_{1j} &=&W_{1j}\left[ \frac{\left( -1\right) ^{r-1}}{\left\| V_{r}\right\|
^{2}...\left\| V_{2}\right\| ^{2}}W_{2}...W_{r^{\prime }}V_{r}...V_{2}\right]
\\
&=&W_{1j}T,
\end{eqnarray*}
where $T$ is defined to be the element $T\in P_{\overline{x}}^{+}M$ shown in
the square brackets above. We note that 
\begin{equation}
T=\frac{-1}{\left\| W_{1j}\right\| ^{2}}W_{1j}V_{1j}\ \text{for each }j,
\label{Telement}
\end{equation}
so that it is an element of $P_{\overline{x}}^{+}M$ of degree at most two.
Choose an orthonormal basis $\left\{ f_{1},...,f_{n}\right\} $ of $T_{%
\overline{x}}M$ such that 
\begin{eqnarray*}
f_{1} &=&\frac{W_{11}}{\left\| W_{11}\right\| }, \\
V_{11} &=&c_{1}f_{1}+c_{2}f_{2}
\end{eqnarray*}
with $c_{1},c_{2}\in \mathbb{R}$. From Equation (\ref{Telement}) with $j=1$,
it is easily seen that 
\begin{equation*}
T=t_{1}+t_{2}f_{1}f_{2}
\end{equation*}
for some constants $t_{1},t_{2}\in \mathbb{R}$.

Suppose that for some $j>1$, $W_{1j}$ is not contained in the space spanned
by $f_{1}$ and $f_{2}$, so that 
\begin{equation*}
W_{1j}=k_{1}f_{1}+k_{2}f_{2}+q,
\end{equation*}
where $q$ is a nonzero vector orthogonal to $f_{1}$ and $f_{2}$. Then 
\begin{eqnarray*}
V_{1j} &=&W_{1j}T \\
&=&t_{1}W_{1j}+t_{2}W_{1j}f_{1}f_{2} \\
&=&\left( t_{1}W_{1j}-t_{2}k_{1}f_{2}+t_{2}k_{2}f_{1}\right)
+t_{2}qf_{1}f_{2},
\end{eqnarray*}
which is a vector if and only if $t_{2}=0$. Thus, if for some $j>1$ there is
a vector $W_{1j}$ orthogonal to the space spanned by $V_{11}$ and $W_{11}$,
then $V_{1k}=tW_{1k}$ for some nonzero real number $t$, for all $k$. Since
the dimension $n$ is even, this implies that the vector fields $V_{1}$ and $%
W_{1}$ have the same index at the origin.

If on the other hand, $W_{1j}$ is contained in the span of $f_{1}$ and $%
f_{2} $ for each $j$, then the dimension $n$ is $2$, and the map 
\begin{equation*}
W_{1j}\longmapsto W_{1j}\left( t_{1}+t_{2}f_{1}f_{2}\right) ,
\end{equation*}
which transforms each $W_{1j}$ to $V_{1j}$, induces an
orientation-preserving linear transformation, whose matrix in the basis $%
\left\{ f_{1},f_{2}\right\} $ is $\left( 
\begin{array}{ll}
t_{1} & -t_{2} \\ 
t_{2} & t_{1}
\end{array}
\right) $. Thus, in this case, the indices of the vector fields $V_{1}$ and $%
W_{1}$ are the same as well.
\end{proof}

\begin{theorem}
\label{PinorHopf}For any section $\beta \in \Gamma \left( P^{-}M\right) $
with nondegenerate zeros, the Euler characteristic satisfies 
\begin{equation*}
\chi \left( M\right) =\sum_{\overline{x}}\mathrm{ind}_{\beta }\left( 
\overline{x}\right) 
\end{equation*}
where the sum is taken over all the zeros $\overline{x}$ of $\beta $.
\end{theorem}

\begin{proof}
Consider the proper perturbation $Z$ of the de Rham operator 
\begin{equation*}
D=d+d^{*}:\Gamma \left( M,\Lambda ^{\mathrm{even}}T^{*}M\otimes \mathbb{C}%
\right) \rightarrow \Gamma \left( M,\Lambda ^{\mathrm{odd}}T^{*}M\otimes %
\mathbb{C}\right) 
\end{equation*}
defined by 
\begin{equation*}
Z_{\beta }=\gamma \otimes c\left( \beta \right) ^{*}\gamma :\mathbb{S}%
\otimes \mathbb{S}^{*}\rightarrow \mathbb{S}\otimes \mathbb{S}^{*},
\end{equation*}
where we have used the isomorphism $\Lambda ^{\bullet }T^{*}M\otimes %
\mathbb{C}\cong \mathbb{S}\otimes \mathbb{S}^{*}$. In a neighborhood $U$ of
a particular zero $\overline{x}$, we write 
\begin{equation*}
\left. \beta \right| _{U}=\left. V_{1}V_{2}...V_{r}\right| _{U},
\end{equation*}
where $V_{1}$ has an isolated nondegenerate zero at $\overline{x}$ and the
vector fields $V_{2},...,V_{r}$ are nonzero at $\overline{x}$. Observe that 
\begin{eqnarray*}
Z_{\beta } &=&r\left( \beta \right) \left( -1\right) ^{p} \\
&=&\left( V_{1}^{\flat }\wedge +i\left( V_{1}\right) \right) \left(
V_{2}^{\flat }\wedge +i\left( V_{2}\right) \right) ...\left( V_{r}^{\flat
}\wedge +i\left( V_{r}\right) \right) 
\end{eqnarray*}
on $p$-forms (see Section \ref{deRham}). We now use Theorem\ref{masterpiece0}
to calculate the index of $D$, which is the Euler characteristic. Without
changing the local index $\mathrm{ind}_{\mathbb{R}^{n}}\left( D\left( 
\overline{x}\right) \right) $, we will deform it in a particularly simple
way near $\overline{x}$. We locally deform the vector fields $V_{2},...,V_{r}
$ smoothly to $\frac{V_{2}}{\left\| V_{2}\right\| }$ (while keeping $\beta $
invertible on $U\diagdown \left\{ \overline{x}\right\} $), and thus the
product $V_{2}...V_{r}$ is deformed to $\left( -1\right) ^{r\left(
r-1\right) /2}$. Then $\beta $ has been deformed to 
\begin{equation*}
\left. \widetilde{\beta }\right| _{U}=\left. \left( -1\right) ^{r\left(
r-1\right) /2}V_{1}\right| _{U}.
\end{equation*}
The perturbation proof of the Poincar\'{e}-Hopf theorem implies that the
difference in dimensions of the kernels of the operators $\left. \left(
D+sZ\right) ^{\pm }\right| _{U}$ for large $s$ is the index of the vector
field $\pm V_{1}$ at $\overline{x}$, which is the same as the index of the $%
V_{1}$ at $\overline{x}$, which is by definition $\mathrm{ind}_{\beta
}\left( \overline{x}\right) $.
\end{proof}

\begin{remark}
The theorem above implies that the odd pinor bundle $\mathrm{Pin}^{-}\left(
TM\right) $ has a section if and only if the Euler characteristic of $M$ is
zero.
\end{remark}

\begin{example}
Consider the standard Dirac operator on the (trivial) spin$^{c}$ bundle $%
\mathbb{S}=\mathbb{C}^{2^{m}}$ over an even-dimensional sphere $%
S^{2m}\subset \mathbb{R}^{2m+1}$. Consider the section $\beta \in \Gamma
\left( P^{-}M\right) $ defined by 
\begin{equation*}
\beta =p\left( E_{1}\right) p\left( E_{2}\right) ...p\left( E_{2m+1}\right) ,
\end{equation*}
where $\left\{ E_{1},...,E_{2m+1}\right\} $ is the standard basis of vector
fields in $\mathbb{R}^{2m+1}$, and where $p_{x}:T_{x}\mathbb{R}
^{2m+1}\rightarrow T_{x}S^{2m}$ is the orthogonal projection. (Note that a
similar bundle map may be constructed on any spin$^{c}$ submanifold of $%
\mathbb{R}^{2m+1}$.) Observe that $\beta \in \Gamma \left( P^{-}M\right) $,
and its zeros are the points of intersection of the axes in $\mathbb{R}
^{2m+1}$ with $S^{2m}$. Note that these zeros are nondegenerate. To see
this, renumber the axes so that the zero in question is the axis parallel to 
$E_{1}$. Near this zero, the vector fields $p\left( E_{2}\right)
,...,p\left( E_{2m+1}\right) $ are nonzero, and $p\left( E_{1}\right) =\pm %
\left( \sin r\right) \partial _{r}$, where $r$ is the geodesic radial
coordinate. Since $p\left( E_{1}\right) $ is locally the gradient of $\mp
\cos \left( r\right) $, which is a Morse function near $r=0$, $p\left(
E_{1}\right) $ has nondegenerate zeros. Also, note that the index of $%
p\left( E_{1}\right) $ is $1$, since it is either a source or sink. Thus $%
\beta $ has nondegenerate zeros. We now compute the index of $\beta $ at
each zero. At one of the zeros of $p\left( E_{j}\right) $, we write 
\begin{eqnarray*}
\beta  &=&p\left( E_{1}\right) p\left( E_{2}\right) ...p\left(
E_{2m+1}\right)  \\
&=&p\left( F_{j}\right) p\left( E_{1}\right) ...p\left( E_{j-1}\right)
p\left( E_{j+1}\right) ...p\left( E_{2m+1}\right)  \\
&=&\pm p\left( F_{j}\right) p\left( \pm E_{1}\right) ...p\left(
E_{j-1}\right) p\left( E_{j+1}\right) ...p\left( E_{2m+1}\right) ,
\end{eqnarray*}
where $F_{j}$ is a vector field with two nondegenerate zeros (at the zeros
of $E_{j}$) whose index is $\left( -1\right) ^{j-1}$ at each of those zeros.
We have used Lemma \ref{CommutingLemma}. The form of $\beta $ above implies
that the index of $\beta $ at each of the two zeros of $p\left( E_{j}\right) 
$ is $\left( -1\right) ^{j-1}$. This verifies Theorem \ref{PinorHopf}, which
implies that 
\begin{eqnarray*}
\chi \left( S^{2m}\right)  &=&\sum_{\overline{x}}\mathrm{ind}_{\beta }\left( 
\overline{x}\right)  \\
&=&\sum_{j=1}^{2m+1}\left[ \left( -1\right) ^{j-1}+\left( -1\right)
^{j-1}\right]  \\
&=&2.
\end{eqnarray*}
\end{example}

\subsection{The induced index on submanifolds\label{submanifoldSection}}

Let $F$ be a real, oriented vector bundle over $M$ with odd rank. Let $%
E=E^{+}\oplus E^{-}$ be a graded, self-adjoint Clifford module over the
total space of $F$. Identify the zero section of $F$ with $M$; the inclusion 
$TM\subset TF$ then induces a graded $\mathbb{C}\mathrm{l}\left( TM\right) $%
-action on $E$ over $M$.

\begin{proposition}
\label{Cobordism}The index of the Dirac operator associated to the $\mathbb{C%
}\mathrm{l}\left( TM\right) $-module $E$ over $M$ is zero.
\end{proposition}

\begin{proof}
Let $\omega $ be a nonzero vertical volume form that induces the orientation
of $F$. Then $Z=ic(\omega ):E^{\pm }\rightarrow E^{\mp }$ is a nonsingular,
self-adjoint bundle map that anticommutes with Clifford multiplication by
sections of $TM$. By Remark \ref{IndexZeroRemark}, the index of the
corresponding Dirac operator is zero.
\end{proof}

Clearly, the same result would apply to an oriented submanifold of odd
codimension in a manifold endowed with a given Clifford module; the vector
bundle $F$ is the normal bundle of the submanifold. It would also apply to a
component of the boundary of a manifold with boundary.

\begin{remark}
This resembles the cobordism invariance of the index of Dirac operators. See 
\cite[Chapter XVII]{Palais}. A perturbation proof of the standard cobordism
invariance result is given by M. Braverman in \cite{Bra5} and \cite{Bra6}.
In our case, the bundle map is an odd endomorphism, and we do not require
that the manifold with boundary be compact.
\end{remark}

\section{Appendix: Example when no localization occurs\label{NoLocalization}}

In this section, we give an example of a perturbation of the Dirac operator
that yields first order terms and where no localization occurs. This
motivates the requirement that $ZD+DZ$ be a bundle map in Section \ref
{perturbDiracOpsSection}.

Let $D_{s}:\Gamma \left( S^{1},\mathbb{C}\right) \rightarrow \Gamma \left(
S^{1},\mathbb{C}\right) $ be defined by 
\begin{equation*}
D_{s}f:=i\frac{df}{d\theta }+s\left( \sin \theta \right) f.
\end{equation*}
Then $D_{s}$ is a family of essentially self-adjoint differential operators.
We solve the equation $D_{s}f=\lambda f$ by separating variables. We
conclude that the eigenvalues of $D_{s}$ are $\lambda _{\pm n}=\pm n$, and
the corresponding orthonormal eigenfunctions are 
\begin{equation*}
f_{\pm n}\left( \theta ,s\right) =\frac{1}{\sqrt{2\pi }}e^{i\left( \mp
n\theta -s\cos \theta \right) }.
\end{equation*}
Clearly, as $s\rightarrow \infty $, eigenvalues $\lambda _{\pm n}$ are
fixed, and the eigenfunctions do not localize in the usual sense since the
magnitude of each eigenfunction stays constant: $\left\vert f_{n}\left(
\theta ,s\right) \right\vert =\frac{1}{\sqrt{2\pi }}$.

In this case, $D=i\frac{d}{d\theta }$, and $Z=\sin \theta $ can be thought
of as Clifford action by the complex vector field $-i\sin \theta \frac{%
\partial }{\partial \theta }$. The anticommutator of $D$ and $Z$ is 
\begin{equation*}
DZ+ZD=i\left( 2\cos \theta \frac{d}{d\theta }-\sin \theta \right) ,
\end{equation*}
and the operator 
\begin{equation*}
D_{s}^{2}=-\frac{d^{2}}{d\theta ^{2}}+s\left( 2i\cos \theta \frac{d}{d\theta 
}-i\sin \theta \right) +s^{2}\cos ^{2}\theta
\end{equation*}
contains a first order term, as expected.

\section{Appendix: Graded Clifford Bundles and Dirac operators\label%
{Preliminary}}

We recall several well known facts about Clifford structures on manifolds of
even and odd dimensions. The books \cite{B-G-V}, \cite{L-M}, and \cite{Roe}
are standard references.

\subsection{The Clifford bundle}

As before, $(M,g)$ is an oriented Riemannian manifold of $\dim M=n$. For any 
$x\in M$, we denote by $\mathrm{Cl}\left( T_{x}M\right) $ the Clifford
algebra of the tangent space $T_{x}M$. The spaces $T_{x}M$ and $T_{x}^{\ast
}M$ are canonically isomorphic, using the chosen Riemannian metric.

The Clifford algebra $\mathrm{Cl}\left( T_{x}M\right) $ is the direct sum of
even and odd components, denoted $\mathrm{Cl}^{+}\left( T_{x}M\right) $ and $%
\mathrm{Cl}^{-}\left( T_{x}M\right) $.

The \textit{Clifford bundle} $\mathrm{Cl}(TM)$ of $M$ is the $Z_{2}$-graded
bundle over $M$ whose fiber at $x\in M$ is $\mathrm{Cl}\left( T_{x}M\right) $
(\cite{B-G-V}, 3.30). We will denote the \textit{complexified Clifford
algebra} $\mathrm{Cl}\left( T_{x}M\right) \otimes \mathbb{C}$ by $\mathbb{C}%
\mathrm{l}\left( T_{x}M\right) $ and the \textit{complexified Clifford bundle%
}\ $\mathrm{Cl}\left( TM\right) \otimes \mathbb{C}$ by $\mathbb{C}\mathrm{l}%
\left( M\right) $.

The Levi-Civita connection $\nabla ^{TM}$ induced by the Riemannian metric $%
g $ extends canonically to a connection on $\mathbb{C}\mathrm{l}\left(
T_{x}M\right) $ compatible with the grading and the Clifford multiplication.

\subsection{Clifford modules\label{CliffordModuleSection}}

A \textit{graded} \textit{self-adjoint} \textit{Clifford module} (\cite
{B-G-V}, 3.32) on a manifold $M$ is a $Z_{2}$-graded complex vector bundle $%
E=E^{+}\oplus E^{-}$ together with a bundle endomorphism 
\begin{equation*}
c:TM\rightarrow \text{\textrm{End}}(E).
\end{equation*}
such that the following properties hold: for any $x\in M$ and any vectors $%
v, $ $w\in T_{x}M$

(i) $c(v):E_{x}^{\pm }\rightarrow E_{x}^{\mp }$ is a graded action;

(ii) $c(v)c(w)+c(w)c(v)=-2g_{x}(v,w)\mathbf{1}$, where $g_{x}$ is the metric
on $T_{x}M$;

(iii) the bundle $E$ is equipped with a Hermitian metric such that the
subbundles $E^{+}$ and $E^{-}$ are orthogonal and the operator $c(v)$ is
skew-adjoint ;

(iv) $E$ is equipped with a grading-preserving Hermitian connection $\nabla
=\nabla ^{E}$ satisfying 
\begin{equation*}
\left[ \nabla _{V}^{E},c(W)\right] =c\left( \nabla _{V}^{TM}W\right) ,
\end{equation*}
for arbitrary vector fields $V$ and $W$ on $M$. This connection is called a 
\textit{Clifford connection }(\cite{B-G-V}, 3.39). Clifford connections
always exist (\cite{B-G-V}, 3.41).

\subsection{Twisted Clifford modules}

Given a Clifford module $E$ and a vector bundle $F$ over $M$, we can
construct the \textit{twisted Clifford module }$E\otimes F$ \textit{obtained
from }$E$ \textit{by twisting with} $F.$ The Clifford action on $E\otimes F$
is given by $c(v)\otimes 1$. Given a connection $\nabla ^{F}$ on $F$ we can
define the \textit{product connection }$\nabla ^{E}\otimes 1+1\otimes \nabla
^{F}$ on $E\otimes F$.

\subsection{The chirality operator and the induced grading on $E$\label%
{chiralitySection}}

Let $e_{1},...,e_{n}$ be an oriented orthonormal basis of $T_{x}M$. We
consider the element 
\begin{equation*}
\gamma =i^{k}c(e_{1})...c(e_{n})\in \mathrm{End}\left( E_{x}\right) \mathbf{,%
}
\end{equation*}
where $k=n/2$ if $n$ is even and $k=\left( n+1\right) /2$ if $n$ is odd.
This element is independent of the choice of basis and anticommutes with any 
$c(v)$ where $v\in T_{x}M$ if $n$ is even and commutes if $n$ is odd.
Moreover, $\gamma ^{2}=\mathbf{1}$ (\cite{B-G-V}, 3.17). The \textit{%
chirality operator }$\gamma $ is the section of $\mathrm{End}\left( E\right) 
$ that restricts to the element above on each fiber. The bundle map $\gamma $
has eigenvalues $\pm 1$, and we can define subbundles 
\begin{equation*}
E_{\gamma }^{\pm }=\left\{ \alpha \in E:\gamma \alpha =\pm \alpha \right\} .
\end{equation*}
The grading $E=E_{\gamma }^{+}\oplus E_{\gamma }^{-}$ is called the\textit{\
grading} \textit{induced by }$\gamma $ on $E$ or the \textit{natural\ grading%
} on $E$.

\subsection{The Dirac operator\label{DiracOpSection}}

The \textit{Dirac operator } $D:\Gamma \left( M,E\right) \rightarrow \Gamma
\left( M,E\right) $ associated to a Clifford connection $\nabla ^{E}$ is
defined by the following composition 
\begin{equation*}
\Gamma \left( M,E\right) \overset{\nabla ^{E}}{\longrightarrow }\Gamma
\left( M,T^{\ast }M\otimes E\right) \overset{c}{\longrightarrow }\Gamma
\left( M,E\right) .
\end{equation*}
In local coordinates this operator may be written as 
\begin{equation*}
D=\sum_{i=1}^{n}c\left( dx_{i}\right) \nabla _{\partial _{i}}:\Gamma \left(
M,E^{\pm }\right) \longrightarrow \Gamma \left( M,E^{\mp }\right) .
\end{equation*}
This is a first order elliptic operator. Moreover, it is formally
self-adjoint and essentially self-adjoint with the initial domain smooth,
compactly supported sections (\cite{B-G-V}, p. 119). Its \textit{principal
symbol }is given by 
\begin{equation*}
\sigma _{D}\left( x,\xi \right) =ic\left( \xi \right) :\Gamma \left(
M,T_{x}^{\ast }M\right) \rightarrow \mathrm{End}\left( E_{x}\right) .
\end{equation*}

\subsection{The $\mathrm{\mathrm{spin}}$ and $\mathrm{\mathrm{spin}}^{c}$
bundles\label{Spinandspinc}}

Let $M$ be an even dimensional oriented manifold with spin structure, and
let $\mathbb{S}=\mathbb{S}^{+}\oplus \mathbb{S}^{-}$ be a complex spinor
bundle over $M$ with the grading induced by $\gamma $. It is a minimal
Clifford module; i.e. for any other Clifford module $E$ over $M$ there is a
vector bundle $F$ such that we have an isomorphism of Clifford modules 
\begin{equation*}
E\cong \mathbb{S}\otimes F,
\end{equation*}
where $F=\mathrm{\mathrm{Hom}}_{\mathbb{C}\mathrm{l}(TM)}\left( \mathbb{S}%
,E\right) $ and the Clifford action is trivial on the second factor (\cite
{B-G-V}, sect. 3.3).

If the dimension of $M$ is odd, then any Clifford module $E$ over $M$ is
isomorphic to 
\begin{equation*}
E\cong \left( \mathbb{S}\otimes F_{1}\right) \oplus \left( \mathbb{S}\otimes
F_{2}\right)
\end{equation*}
where $v\in \Gamma (M,TM)$ acts on $\mathbb{S}\otimes F_{1}$ by $c\left(
v\right) \otimes \mathbf{1}$ and on $\mathbb{S}\otimes F_{2}$ by $c\left(
-v\right) \otimes \mathbf{1}$. In the odd case we denote the first action $%
c^{+}\left( v\right) \otimes \mathbf{1}$ and the second action $c^{-}\left(
v\right) \otimes \mathbf{1}$.

A connection $\nabla ^{E}$ on the twisted Clifford module $E=\mathbb{S}%
\otimes F$ is a Clifford connection if and only if 
\begin{equation*}
\nabla ^{E}=\nabla ^{\mathbb{S}}\otimes 1+1\otimes \nabla ^{F}
\end{equation*}
for some connection $\nabla ^{F}$ on $F$.

Note that there are global obstructions to the existence of complex spinor
bundles (see \cite{B-G-V}, 3.34); however, locally the decompositions above
always exist.

Every $\mathrm{\mathrm{spin}}$ manifold and every almost complex manifold
has a canonical $\mathrm{\mathrm{spin}}^{c}$ structure. In addition, every
oriented, compact manifold of dimension $\leq 3$ is $\mathrm{\mathrm{spin}}%
^{c}$ \cite{Nic}\textbf{.}

\subsection{Local classification of gradings}

The above results apply to Dirac operators over bundles with the natural
grading --- that induced directly from the grading on complex spinors. The
following lemma classifies all possible gradings for Clifford
representations (and thus Dirac operators).

Let $V$ be an oriented, real Euclidean vector space, and let $E=E^{+}\oplus
E^{-}$ be a complex vector space that is a graded $\mathbb{C}\mathrm{l}
\left( V\right) $-module. Let $\mathbb{S}$ denote\ the irreducible
representation space of $\mathbb{C}\mathrm{l}\left( V\right) $. Let $c\left(
v\right) $ denote the Clifford multiplication by $v\in V$on $\mathbb{S}= %
\mathbb{S}^{+}\oplus \mathbb{S}^{-}$ when $V$ is even-dimensional, and let $%
c^{+}\left( v\right) :\mathbb{S}\rightarrow \mathbb{S}$ and $c^{-}\left(
v\right) :=-c^{+}\left( v\right) :\mathbb{S}\rightarrow \mathbb{S}$ be the
Clifford multiplications that generate the two nonequivalent irreducible
representations of $\mathbb{C}\mathrm{l}\left( V\right) $ when $V$ is
odd-dimensional. We let 
\begin{equation*}
c^{E}\left( v\right) :E^{\pm }\rightarrow E^{\mp }
\end{equation*}
denote the Clifford action by $v$ on the graded vector space $E$.

\begin{lemma}
\label{GradingClassification} With the notation described above, there
exists a complex vector space $W$ such that

\begin{enumerate}
\item  $E\cong \mathbb{S}\otimes W$, where $W\cong \mathrm{Hom}_{\mathbb{C}%
\mathrm{l}\left( V\right) }\left( \mathbb{S},E\right) $.

\item  If the dimension of $V$ is even, then the associated Clifford action
on $\mathbb{S}\otimes W$ is $c\left( v\right) \otimes \mathbf{1}$. In
addition, $W=W^{+}\oplus W^{-}$ for some orthogonal subspaces $W^{\pm }$ of $%
W$, and $E^{\pm }\cong \left( \mathbb{S}^{+}\otimes W^{\pm }\right) \oplus
\left( \mathbb{S}^{-}\otimes W^{\mp }\right) $.\label{even case}

\item  If the dimension of $V$ is odd\label{odd case}, then there exists an
orthogonal decomposition $W=W^{\prime }\oplus W^{\prime }$, such that 
\begin{equation*}
E\cong \mathbb{S}\otimes W\cong \left( \mathbb{S}\otimes W^{\prime }\right)
\oplus \left( \mathbb{S}\otimes W^{\prime }\right) ,
\end{equation*}
the induced action of $c^{E}\left( v\right) $ on $\left( \mathbb{S}\otimes
W^{\prime }\right) \oplus \left( \mathbb{S}\otimes W^{\prime }\right) $ is
given by 
\begin{equation*}
\left( c^{+}\left( v\right) \otimes \mathbf{1,\,}c^{-}\left( v\right)
\otimes \mathbf{1}\right) =\left( c^{+}\left( v\right) \otimes \mathbf{1,\,}%
c^{+}\left( v\right) \otimes -\mathbf{1}\right) ,
\end{equation*}
and 
\begin{equation*}
E^{\pm }\cong \mathbb{S}\otimes \mathrm{span}\left\{ \left( w,\pm w\right)
\in W^{\prime }\oplus W^{\prime }\,|\,w\in W^{\prime }\right\} .
\end{equation*}
\end{enumerate}
\end{lemma}

\begin{proof}
The first statement follows directly from the general facts about the
representation theory of Clifford algebras; see Section \ref{Spinandspinc}.%
\newline
To prove (\ref{even case}), observe that the even part of the Clifford
algebra $\mathbb{C}\mathrm{l}^{+}\left( V\right) $ acts by endomorphisms on $
E^{+}$ and $E^{-}$. This leads to a representation of $\mathrm{spin}
^{c}\left( n\right) \subset \mathbb{C}\mathrm{l}^{+}$ on $E^{+}$ and $E^{-}$
. There are exactly two nonequivalent irreducible representations of $%
\mathrm{spin}^{c}\left( n\right) $, given by the actions of $\mathbb{C}%
\mathrm{l}^{+}$ on $\mathbb{S}^{+}$ and $\mathbb{S}^{-}$; see \cite[p. 432]
{Nic}. Thus, there are complex vector spaces $W^{+}$ and $W^{-}$ such that $%
E^{+}\cong \left( \mathbb{S}^{+}\otimes W^{+}\right) \oplus \left( \mathbb{S}
^{-}\otimes W^{-}\right) $, which implies that $E^{-}\cong \left( \mathbb{S}%
^{-}\otimes W^{+}\right) \oplus \left( \mathbb{S}^{+}\otimes W^{-}\right) $.

To prove (\ref{odd case}), where $V$ is odd-dimensional, observe that there
is a unique irreducible representation of $\mathrm{spin}^{c}\left( n\right) $
, given by the action of $\mathbb{C}\mathrm{l}^{+}$ on $\mathbb{S}$. We have 
$E=E^{+}\oplus E^{-}\cong \left( \mathbb{S}\otimes W^{+}\right) \oplus
\left( \mathbb{S}\otimes W^{-}\right) \cong \mathbb{S}\otimes \left(
W^{+}\oplus W^{-}\right) $ with the Clifford action on the last term being $%
c^{E}(v)=c^{+}(v)\otimes J$. Here $J:W^{\pm }\rightarrow W^{\mp }$. The
operator $J$ must be Hermitian and squares to identity (since $c^{E}(v)$ is
skew-hermitian and squares to $-\mathbf{1}$).

Now let $\left\{ e_{1},...,e_{k}\right\} $ be an orthonormal basis of $W^{+}$
; then $\left\{ Je_{1},...,Je_{k}\right\} $ must be an orthonormal basis of $%
W^{-}.$ Thus $W$ has an orthonormal basis $\left\{
e_{1},...,e_{k},Je_{1},...,Je_{k}\right\} .$ We introduce a new
decomposition of $W=W_{1}\oplus W_{2}$, where $W_{1}=\mathrm{span}\left\{
e_{1}+Je_{1,...,}e_{k}+Je_{k}\right\} $ and $W_{2}=\mathrm{span}\left\{
e_{1}-Je_{1},...,e_{k}-Je_{k}\right\} .$ Then $J$ is the identity on $W_{1}$
and minus the identity on $W_{2}$, and we can decompose $E$ as 
\begin{equation*}
E\cong \left( \mathbb{S}\otimes W_{1}\right) \oplus \left( \mathbb{S}\otimes
W_{2}\right) ,
\end{equation*}
where $c^{E}(v)$ acts by $\left( c^{+}(v)\otimes \mathbf{1},c^{+}(v)\otimes -%
\mathbf{1}\right) $. The conclusion (\ref{odd case}) follows from the
observation that $W^{\prime }:=W_{1}\cong W_{2}$, where the isomorphism maps
each $e_{m}+Je_{m}$ to $e_{m}-Je_{m}$.
\end{proof}

\subsection{Global classification of gradings\label%
{globalClassificationSection}}

Let $E=E^{+}\oplus E^{-}$ be a graded, self-adjoint, Hermitian $\mathbb{C}%
\mathrm{l}\left( TM\right) $ -module over $M$. Suppose that $M$ is spin$^{c}$%
. Choose a spin$^{c}$ structure on $M$. Let $\mathbb{S}$ be the
corresponding spin$^{c}$ bundle over $M$, so that the representation of $%
\mathbb{C}\mathrm{l}\left( TM\right) $ is irreducible on $\mathbb{S}$ and $%
\mathbb{C}\mathrm{l}\left( TM\right) \cong \mathrm{End}\left( \mathbb{S}%
\right) $. Let $c\left( v\right) $ denote the Clifford multiplication by $%
v\in TM\otimes \mathbb{C}$ on $\mathbb{S}$; $\mathbb{S}=\mathbb{S}^{+}\oplus %
\mathbb{S}^{-}$ if $M$ is even-dimensional. Let $c^{+}\left( v\right) :%
\mathbb{S}\rightarrow {}\mathbb{\ S}$ and $c^{-}\left( v\right)
:=-c^{+}\left( v\right) :\mathbb{S}\rightarrow \mathbb{S}$ be the Clifford
multiplications that generate the two irreducible representations of $%
\mathbb{C}\mathrm{l}\left( TM\right) $ when $M$ is odd-dimensional. We let 
\begin{equation*}
c^{E}\left( v\right) :E^{\pm }\rightarrow E^{\mp }
\end{equation*}
denote the Clifford multiplication by $v$ on the graded vector bundle $E$.

\begin{corollary}
(Classification of bundles of graded Clifford modules)\label%
{GradedBundleClassification} If $M$ is spin$^{c}$ with the above notation,
there exists a complex vector bundle $W$ such that

\begin{enumerate}
\item  $E\cong \mathbb{S}\otimes W$, where $W\cong \mathrm{Hom}_{\mathbb{C}%
\mathrm{l}\left( M\right) }\left( \mathbb{S},E\right) $.

\item  When $M$ is even-dimensional, the associated Clifford action on $%
\mathbb{S}\otimes W$ is $c\left( v\right) \otimes \mathbf{1}$. Then $%
W=W^{+}\oplus W^{-}$ for some orthogonal vector subbundles $W^{\pm }$ of $W$
, and $E^{\pm }\cong \left( \mathbb{S}^{+}\otimes W^{\pm }\right) \oplus
\left( \mathbb{S}^{-}\otimes W^{\mp }\right) $.

\item  When $M$ is odd-dimensional, there exists an orthogonal decomposition 
$W=W^{\prime }\oplus W^{\prime }$, such that 
\begin{equation*}
E\cong \mathbb{S}\otimes W\cong \left( \mathbb{S}\otimes W^{\prime }\right)
\oplus \left( \mathbb{S}\otimes W^{\prime }\right) ,
\end{equation*}
the induced action of $c^{E}\left( v\right) $ on $\left( \mathbb{S}\otimes
W^{\prime }\right) \oplus \left( \mathbb{S}\otimes W^{\prime }\right) $ is
given by 
\begin{equation*}
\left( c^{+}\left( v\right) \otimes \mathbf{1,\,}c^{-}\left( v\right)
\otimes \mathbf{1}\right) =\left( c^{+}\left( v\right) \otimes \mathbf{1,\,}%
c^{+}\left( v\right) \otimes -\mathbf{1}\right) ,
\end{equation*}
and 
\begin{equation*}
E^{\pm }\cong \mathbb{S}\otimes \mathrm{span}\left\{ \left( w,\pm w\right)
\in W^{\prime }\oplus W^{\prime }\,|\,w\in W^{\prime }\right\} .
\end{equation*}
\newline
\end{enumerate}

Finally, if $M$ is not spin$^{c}$, then the relevant facts above are true
locally but not globally; that is, the bundles $\mathbb{S}$, $W$, $W^{\pm }$%
, $W^{\prime }$ can be defined on a sufficiently small neighborhood of any
given point of $M$, and the properties above are true over that
neighborhood, but $\mathbb{S}$ cannot be extended to a globally defined spin$%
^{c}$ bundle.
\end{corollary}

\begin{proof}
The fact that $E\cong \mathbb{S}\otimes W$ in both the odd and even cases
follows by setting $W=\mathrm{Hom}_{\mathbb{C}\mathrm{l}}\left( \mathbb{S}%
,E\right) $, the bundle maps from $\mathbb{S}$ to $E$ that are $\mathbb{C}%
\mathrm{l}\left( TM\right) $-equivariant. (In the odd case, one must fix an
irreducible representation $c^{+}$ on $\mathbb{S}$ once and for all.) The
isomorphism $\mathbb{\ S}\otimes W\rightarrow E$ is given by $s\otimes
w\mapsto w\left( s\right) $.

In the even case, observe that the bundles $W^{\pm }$ may be defined
globally by noting that for example $\left( \mathbb{S}^{+}\otimes W^{\pm
}\right) =E^{\pm }\cap \left( \mathbb{S}^{+}\otimes W\right) $, where we
have abused notation using the isomorphism $\mathbb{S}\otimes W\rightarrow E$%
.

In the odd case, observe that in the proof of Lemma \ref
{GradingClassification}, the representation theory alone determines the
bundles $W_{1}$ and $W_{2}$ from $W$, and the construction of $W^{\prime }$
is canonical. The result follows.
\end{proof}

\end{document}